%% file: laminations_hyperfinies.tex
\documentclass[10pt]{article}

\usepackage{amsthm,amsmath,amsfonts,amssymb}
\usepackage[latin1]{inputenc}
\usepackage{epsfig,color}
\usepackage{float}

\usepackage{pb-diagram}
\usepackage[all]{xy}

\newtheorem{thm}{Théorème}[section]
\newtheorem{cor}[thm]{Corollaire}
\newtheorem{lem}[thm]{Lemme}
\newtheorem{prop}[thm]{Proposition}

\theoremstyle{definition}\newtheorem{defn}[thm]{Définition}
\theoremstyle{definition}\newtheorem{ex}[thm]{Exemple}
\theoremstyle{definition}\newtheorem{rem}[thm]{Remarque}

\newcommand{\R}{\mathbb R}
\newcommand{\C}{\mathbb C}
\newcommand{\Z}{\mathbb Z}
\newcommand{\N}{\mathbb N}

\newcommand{\K}{\mathcal K}
\newcommand{\preuve}{\noindent{\it Démonstration}. }

\newcommand{\mois}{%
\ifcase\month\or
 Janvier\or Février\or Mars\or Avril\or Mai\or Juin\or
 Juillet\or Août\or Septembre\or Octobre\or Novembre\or
Décembre\fi
}
\renewcommand{\today}{\number\day\ \mois\ \number\year}

\title{{\bf Laminations hyperfinies et revêtements}}
\date{\today}

\author{
{\bf Miguel Bermúdez\footnote{Ce travail a été partiellement financé
    par le projet  BFM2001-3280 Ministerio de Ciencia y Tecnología
    (Espagne).}}\\ 
\and
{\bf Gilbert  Hector}\\
}

\begin{document}

\maketitle \renewcommand{\abstractname}{Resumé}

\begin{abstract}
  Nous introduisons la catégorie des espaces boréliens-topologiques,
  appelée BT-catégorie: c'est un cadre naturel pour la classification
  mesurable des feuilletages et des laminations usuels. Nous nous
  intéressons au cas des laminations de dimension deux. Nous prouvons
  les deux résultats principaux suivants:

  \medskip\noindent (1) Une lamination borélienne par plans est la
  suspension d'une action de $\Z^2$ sur un espace de Borel standard si
  et seulement si elle est hyperfinie.

\noindent
(2) Toute lamination par surfaces admettant une structure complexe
parabolique mesurable est moyennable.

\medskip Le troisième résultat est une amélioration de (2) dans le cas
des laminations munies d'une mesure transverse quasi-invariante $\mu$.
L'énoncé est le suivant:

\medskip\noindent (3) Une lamination par plans, cylindres et tores est
$\mu$-moyennable si et seulement si elle admet une métrique de Riemann
plate le long de $\mu$-presque chaque feuille.
\end{abstract}

\renewcommand{\abstractname}{Abstract}
\begin{abstract}
  We introduce the so-called BT-category of borelian-topo\-logical
  spaces: it will be a natural frame for a measurable classification
  of usual foliations and laminations. We focus on the two-dimensional
  case: borelian laminations by surfaces. We prove two main results:

  \medskip\noindent (1) Any borelian lamination by planes is the
  suspension of a $\Z^2$-action on a Borel space iff this lamination
  is hyperfinite.

\noindent
(2) Any borelian lamination by surfaces is amenable if parabolic, i.e.
if it admits a complex structure parabolic on each leaf.

\medskip The third result is an improvement of (2) in case of
laminations endowed with a transverse quasi-invariant measure $\mu$.
The statement is the following:

\medskip\noindent (3) Any borelian lamination by planes, cylinders and
tori is $\mu$-amenable if and only if it admits a metric which is flat
on $\mu$-almost all leafs.
\end{abstract}

\section*{Introduction}
L'objet du présent travail est une étude des laminations par surfaces
d'un point de vue mesurable. Concrètement, étant donnée une lamination
$(X,\mathcal F)$ on considère la double structure donnée par:
\begin{enumerate}
\item[a)] la structure borélienne $\mathcal B$ engendrée par la
  topologie initiale de $X$;

\item[b)] la topologie fine de la lamination, encore notée $\mathcal
  F$.
\end{enumerate}
Ce couple définit une {\em lamination borélienne}; c'est un objet
d'une catégorie que nous désignerons par BT (borélienne-topologique).
Les laminations boréliennes sont alors une généralisation naturelle
des laminations et feuilletages.  Nous avons donc dans ce cadre les
notions naturelles de:

\begin{itemize}
\item {\em Métrique riemannienne et structure conforme:} une métrique
  riemannienne $g$ sur $\mathcal F$ est une famille de métriques sur
  les feuilles de $\mathcal F$ variant de façon borélienne. Toute
  métrique de Riemann détermine une famille borélienne de structures
  conformes sur les feuilles de $(X,\mathcal F)$.

\item {\em Mesure de Lebesgue:} Une métrique de Riemann sur
  $(X,\mathcal F)$ déter\-mine une famille borélienne $\lambda$ de
  mesures sur les feuilles de $\mathcal F$. Cette famille sera appelée
  la {\em mesure de Lebesgue} de $\mathcal F$. La classe d'équivalence
  de $\lambda$, qui est indépendante de la métrique de Riemann
  choisie, sera appelée la {\em classe de Lebesgue}.

\item {\em Triangulation:} ce sera également une famille borélienne de
  triangulations des feuilles de $\mathcal F$. Le sens du mot
  "borélien" dans ce cas sera précisé au \S\ref{sec:tri}.
\end{itemize}

Certains procédés usuels de construction se transposent au cadre des
laminations boréliennes:

\begin{itemize}
\item {\em Suspension:} étant donnée une surface $S$, on construit la
  suspension d'une action borélienne $\phi$ de $\pi_1(S)$ sur un
  espace de Borel standard $T$ de façon analogue au cas continu.
  L'espace total $X$ de cette suspension est un espace de Borel
  standard muni d'une lamination borélienne $\mathcal F$ obtenue par
  passage au quotient de la lamination triviale sur $\tilde S\times
  T$, où $\tilde S$ désigne le revêtement universel de $S$. La
  projection naturelle: $$
  p:X\to S $$
  est une application borélienne,
  continue le long des feuilles, i.e. une application BT. C'est en
  plus un revêtement de $S$ en restriction à chaque feuille. La
  lamination borélienne $(X,\mathcal F)$ est ce que nous appellerons
  un BT-revêtement de $S$.

\item{\em Action d'un groupe de Lie $G$:} une lamination borélienne
  $(X,\mathcal F)$ est définie par une action de $G$ si les feuilles
  de $\mathcal F$ sont les orbites d'une action $$
  \Phi:G\times X\to X
  $$
  supposée BT en un sens évident. La BT-action $\Phi$ induit sur
  chaque feuille de $\mathcal F$ une action continue transitive de
  $G$. On dira que $\Phi$ est localement libre s'il en est ainsi pour
  l'action induite par $\Phi$ sur chaque feuille.
\end{itemize}

Deux dernières notions jouent un rôle fondamental dans la théorie de
systèmes dynamiques: l'hyperfinitude et la moyennabilité. Une
lamination borélienne $(X,\mathcal F)$ sera dite:
\begin{itemize}
\item[i)] {\em Hyperfinie:} si pour toute transversale borélienne $T$,
  la relation borélienne discrète $\mathcal F_T$ induite par $\mathcal
  F$ sur $T$ est hyperfinie. Ceci signifie qu'il existe une suite
  croissante de relations d'équivalence boréliennes finies $\mathcal
  R_n$ sur $T$ telle que $\mathcal F_T=\cup_n\mathcal R_n$.

\item[ii)] {\em Moyennable:} Soit $\lambda$ la classe de Lebesgue de
  $(X,\mathcal F)$ et $L^\infty(\mathcal F,\lambda)$ l'espace des
  fonctions boréliennes $f:X\to \mathcal \R$ qui sont
  $\lambda$-essentielle\-ment bornées sur chaque feuille.  Alors
  $(X,\mathcal F)$ est dite {\em moyennable} s'il existe une
  application linéaire positive $$
  m:L^\infty(\mathcal F,\lambda)\to
  L^\infty(\mathcal F,\lambda) $$
  qui à tout élément $f\in
  L^\infty(\mathcal F,\lambda)$ associe une fonction $m(f)$ basique
  i.e. constante le long des feuilles, et telle que $m(\bf 1)=\bf 1$.
\end{itemize}

Les notions d'hyperfinitude et moyennabilité peuvent être
"relativisées" par rapport à une mesure transverse quasi-invariante
$\mu$ (cf. \S\ref{sec:mes}). En général, on dira qu'une lamination
borélienne $(X,\mathcal F)$ vérifie une propriété mod$(\mu)$ s'il
existe un borélien saturé $X_1\subset X$ de mesure totale tel que la
lamination borélienne induite $(X_1,\mathcal F_1)$ vérifie cette
propriété. On peut parler alors de lamination hyperfinie mod$(\mu)$.
On pourrait parler également de lamination moyennable mod$(\mu)$, mais
nous aurons besoin d'une notion plus subtile, connue dans la
litérature sous le nom de $\mu$-moyennabilité (cf. \cite{AR}).On note
$L^\infty(\mathcal F,\lambda,\mu)$ l'espace des fonctions boréliennes
qui sont $\lambda$-essentielle\-ment bornées sur $\mu$-presque chaque
feuille. On dira que $(X,\mathcal F)$ est $\mu$-moyennable s'il existe
une application linéaire positive $$
m:L^\infty(\mathcal
F,\lambda,\mu)\to L^\infty(\mathcal F,\lambda,\mu) $$
qui assigne à
toute fonction une fonction basique, et telle que $m({\bf 1})=\bf 1$.
La moyennabilité mod$(\mu)$ entraîne clairement la
$\mu$-moyennabilité, mais la réciproque n'est pas évidente. En
revanche il est connu que l'hyperfinitude mod$(\mu)$ est équivalente à
la $\mu$-moyennabilité (cf.  \cite{CFW}). On aimerait un résultat
analogue dans le cadre purement borélien. A ce propos Weiss
conjecturait dans \cite{W} que toute action borélienne d'un groupe
discret moyennable engendrait une relation d'équivalence hyperfinie.
Cette conjecture a été prouvée dans \cite{JKL} pour un groupe de type
fini à croissance polynomiale.

\medskip Le résultat central de ce travail est le suivant. Il peut
être vu comme une extension au cadre des laminations par surfaces du
théorème de Ambrose-Kakutani pour les flots \cite{AK} (voir aussi
\cite{vN}).

\medskip
\begin{thm}\label{thm:hyp-rev-act}
  Soit $(X,\mathcal F)$ une lamination borélienne par surfaces. Les
  trois conditions suivantes sont équivalentes:
  \begin{enumerate}
  \item $(X,\mathcal F)$ est hyperfinie par plans, tores et cylindres;

  \item $(X,\mathcal F)$ est un BT-revêtement du tore $\mathbb T^2$;

  \item $(X,\mathcal F)$ est définie par une BT-action localement
    libre de $\R^2$.
  \end{enumerate}
\end{thm}

\medskip Il est clair que tout BT-revêtement du tore est défini par
une BT-action localement libre de $\R^2$. En effet, on peut relever
sur les feuilles de la lamination l'action naturelle de $\R^2$ sur le
tore. L'implication (3)$\Rightarrow$(1) découle immédiatement des
résultats de \cite{JKL}. Il ne nous reste donc qu'à montrer
l'implication (1)$\Rightarrow$(2). Pour ce faire nous remplaçons
l'hyperfinitude par une condition qui dans ce cadre est équiva\-lente:
l'{\em hypercompacité}. Cette condition nous permettra de construire
un revêtement du tore par approximations successives sur des morceaux
compacts dans les feuilles. Nous traiterons de manière indépendante
les laminations par plans, par cylindres et par tores. Ceci est
possible car les ensembles saturés formés par les feuilles ayant $0$,
$1$ ou $2$ bouts respectivement sont des boréliens (voir \cite{Gh1}).
Nous traiterons en détail le cas des laminations par plans (théorème
\ref{thm:rev}) et laisserons au lecteur le cas des tores et des
cylindres. Ce dernier est particulièrement intéressant du fait que
toute lamination borélienne dont toutes les feuilles ont deux bouts
est automatiquement hyperfinie car c'est l'{\em extension compacte
  d'un flot} \cite{Gh1} (voir aussi le théorème 3.28 de \cite{Bla}).
Enfin le cas des tores est presque trivial.

\medskip Rappelons qu'une surface de Riemann est dite {\em
  parabolique} si son revêtement universel est conformément équivalent
au plan euclidien; ceci implique en particulier que la surface est
homéomorphe à un tore, un cylindre ou un plan. On peut alors se
demander si toute lamination par tores, plans et cylindres admet une
métrique parabolique. Le résultat suivant répond négativement à cette
question.

\begin{thm}\label{thm:lam-moy}
  Soit $(X,\mathcal F)$ une lamination borélienne par surfaces. Si
  $(X,\mathcal F)$ admet une métrique de Riemann parabolique, alors
  elle est moyennable.
\end{thm}

Toute relation d'équivalence hyperfinie est moyennable, mais la
réciproque est au stade de conjecture au moment où nous écrivons ce
texte. En particulier, nous ignorons toujours si toute lamination
parabolique est hyperfinie. Nous sommes pour cette raison obligés
d'introduire une mesure transverse quasi-invariante dans l'énoncé du
théorème suivant, qui découle des théorèmes \ref{thm:lam-moy} et
\ref{thm:hyp-rev-act} et du théorème de Connes-Feldman-Weiss
\cite{CFW}.

\begin{thm}\label{thm:eq-mes}
  Soit $(X,\mathcal F)$ une lamination borélienne orientable par
  surfaces munie d'une mesure transverse quasi-invariante $\mu$. Alors
  les six conditions suivantes sont équivalentes mod$(\mu)$:
  \begin{enumerate}
  \item $(X,\mathcal F)$ est hyperfinie par plans, tores et cylindres;

  \item $(X,\mathcal F)$ est un BT-revêtement du tore;

  \item $(X,\mathcal F)$ est définie par une BT-action localement
    libre de $\R^2$;

  \item $(X,\mathcal F)$ admet une métrique de Riemann plate complète;

  \item $(X,\mathcal F)$ admet une métrique de Riemann parabolique;

  \item $(X,\mathcal F)$ est moyennable par plans, tores et cylindres.
  \end{enumerate}
\end{thm}

Les implications
(1)$\Rightarrow$(2)$\Rightarrow$(3)$\Rightarrow$(4)$\Rightarrow$(5)$\Rightarrow$(6)
sont vraies dans le cas purement borélien, i.e. sans mesure
quasi-invariante. En effet les deux premières découlent du théorème
\ref{thm:hyp-rev-act} et toute BT-action localement libre de $\R^2$
détermine naturellement une métrique de Riemann complète plate sur les
feuilles de $\mathcal F$ invariante par l'action. Une surface de
Riemann plate {\em complète} est toujours parabolique car son
revêtement universel est le plan euclidien, et enfin l'implication
(5)$\Rightarrow$(6) découle du théorème \ref{thm:lam-moy}. En ce qui
concerne l'implication (6)$\Rightarrow$(1), on a déjà souligné que la
moyennabilité entraîne la $\mu$-moyen\-nabilité, et celle-ci
l'hyperfinitude mod($\mu$) d'après \cite{CFW}.

\medskip Pour illustrer la portée des résultats précédents on peut les
rapprocher des théorèmes de classification des actions localement
libres de $\R^2$ sur des variétés fermées de dimension $3$ (cf.
\cite{Ros}). Ceux-ci montrent qu'une variété fermée de dimension $3$
admettant une telle action est un fibré en tores au dessus de $\mathbb
S^1$. Le théorème \ref{thm:hyp-rev-act} a comme conséquence que toutes
ces actions sont BT-isomorphes à un BT-revêtement du tore $\mathbb
T^2$, et sont donc définies par une BT-action de $\R^2$ sur $\mathbb
T^3$. En adaptant les arguments de \cite{AOW} on pourrait sans doute
montrer que ladite action est presque partout continue (relativement à
la mesure de Lebesgue).

\medskip On remarque qu'il existe des feuilletages par plans de
codimension $3$ et classe $C^\infty$ sur des fibrés principaux
orthogonaux qui ne sont pas moyennables pour la mesure de Lebesgue et
ne peuvent donc pas être définis par une BT-action de $\R^2$ (voir
\S\ref{sec:hyp}).  Cet exemple est aussi à rapprocher du théorème de
Rosenberg sur les feuilletages par plans des variétés fermées de
dimension $3$: ceux-ci sont définis par une action continue de $\R^2$
sur le tore $\mathbb T^3$. Ils sont donc bien hyperfinis et
moyennables.

\medskip De la même façon on peut appliquer le théorème
\ref{thm:eq-mes} aux actions localement libres du groupe affine de la
droite réelle sur une variété compacte.  Puisque le groupe affine est
moyennable, le feuilletage est moyennable et donc hyperfini mod$(\mu)$
pour toute mesure transverse quasi-invariante $\mu$ \cite{CFW}. Par le
théorème \ref{thm:eq-mes} ce feuilletage est BT-isomorphe mod$(\mu)$ à
une BT-action de $\R^2$. En particulier il admet une BT-métrique plate
complète, alors que toutes ses métriques continues sont hyperboliques.

\medskip La question de l'existence de métriques complètes plates sur
les feuilletages a été étudiée par divers auteurs dont Candel
\cite{Ca}, Ghys \cite{Gh2} et Glutsyuk \cite{Glu}. Dans \cite{Gh2}
Ghys construit un exemple de feuilletage par surfaces sur une variété
$M$ de dimension $6$ muni d'une métrique continue parabolique $g$
telle que pour toute BT-fonction $u:M\to \R$, la métrique parabolique
$\exp(u)g$ est non plate. Ce feuilletage est moyennable d'après le
théorème \ref{thm:lam-moy}, et par conséquent hyperfini pour la mesure
de Lebesgue. Il admet en vertu du théorème \ref{thm:eq-mes} une
métrique de Riemann mesurable plate sur presque chaque feuille. Mais
cette métrique n'est pas conformément équivalente à la métrique
parabolique d'origine.

\paragraph*{Remerciements.} Le premier auteur remercie Fernando
Alcalde Cuesta pour ses nombreuses remarques et suggestions et son
aide dans la mise au point des différentes définitions de
moyennabilité.

\bigskip
\section{Laminations par surfaces}
\subsection{Espaces de Borel standard}
Un {\em espace mesurable} est un ensemble $X$ muni d'une
$\sigma$-algèbre $\mathcal B$ de sous-ensembles dits mesurables. On
appellera $\mathcal B$ une {\em structure mesurable} sur $X$. Une
fonction $$
f:(X_1,\mathcal B_1)\to (X_2,\mathcal B_2) $$
entre deux
espaces mesurables est dite {\em mesurable} si, pour tout
$A\in\mathcal B_2$ on a $f^{-1}(A)\in \mathcal B_1$.  Une fonction
mesurable ne préserve pas nécessaire\-ment les structures mesurables.
Autrement dit, l'image directe d'un ensemble mesurable n'est pas
forcément mesurable. Une fonction mesurable $f$ sera dite {\em
  bi-mesurable} si elle envoie ensembles mesurables sur ensembles
mesurables. Si $f$ est bijective, alors c'est un {\em isomorphisme
  mesurable}.

Un espace Polonais est un espace topologique homéomorphe à un espace
métrique complet séparable. Un {\em espace de Borel standard} est un
espace mesurable mesurablement isomorphe à un espace Polonais. Les
ensembles mesurables d'un espace de Borel seront appelés {\em
  boréliens}. Une application mesurable entre deux espaces Borel
standard sera dite une {\em application borélienne}. Un espace de
Borel standard infini est mesurablement isomorphe à l'un des deux
espaces de Borel suivants:
\begin{itemize}
\item L'intervalle $[0,1]$;
\item L'ensemble des nombres entiers positifs $\N$.
\end{itemize}

Le résultat suivant a été démontré par Kallman dans \cite{Ka}. Ce sera
un outil essentiel pour notre étude.

\medskip
\begin{thm}[\cite{Ka}] \label{thm:kal-for} Soient $X$ un espace de
  Borel standard et $B$ un espace polonais. Soit $p_X:X\times B\to X$
  la projection sur le premier facteur. Si $A$ est un borélien de
  $X\times B$ et les ensembles $A_x=\{y\in B~|~ (x,y)\in A\}$ $(x\in
  X)$ sont réunion dénombrable de compacts de $B$, alors $p_X(A)$ est
  un borélien de $X$.  En plus, il existe une section borélienne
  $s:p_X(A)\to A$ de $p_X$.
\end{thm}

\medskip Nous avons le corollaire important suivant: \medskip
\begin{thm}\label{thm:kal}
  Soit $X$ et $Y$ deux espaces de Borel standard et soit $f:X\to Y$
  une application borélienne à fibres dénombrables. Alors $f(X)$ est
  un borélien de $Y$ et il existe une section borélienne $s:f(X)\to X$
  de $f$. En particulier, si $f$ est injective, c'est un isomorphisme
  de Borel entre $X$ et $Y$.
\end{thm}
\preuve Le graphe de $f$, i.e. l'ensemble $A=\{(x,y)\in X\times Y~|~
f(x)=y\}$, est un borélien de $X\times Y$ (voir \cite[page 143]{Hal}).
Comme $X$ est mesurablement isomorphe à un espace polonais ($[0,1]$ ou
$\N$), on peut appliquer le théorème précédent à $A$ et à la
projection de $X\times Y$ sur $Y$.  \qed

\subsection{Relations d'équivalence boréliennes discrètes}\label{sec:rel-mes}
Soit $T$ un espace de Borel standard et $\mathcal R$ une relation
d'équivalence définie sur $T$. On dira que $\mathcal R$ est une {\em
  relation d'équivalence borélienne} sur $T$ si le graphe de $\mathcal
R$ (que nous noterons encore $\mathcal R$) est un borélien de $T\times
T$. Une relation d'équivalence borélienne est dite {\em discrète} si
toutes ses classes sont dénombrables. Elle est {\em finie} si toute
ses classes sont finies. On notera $\mathcal R(x)$ la classe
d'équivalence de $x\in T$. Pour un ensemble $A\subset T$ on notera
$sat_{\mathcal R}(A)$ le $\mathcal R$-saturé $A$, i.e. la réunion des
classes d'équivalence de $\mathcal R$ qui rencontrent $A$. En
utilisant le fait que $\mathcal R$ est un borélien de $T\times T$, on
vérifie que le $\mathcal R$-saturé de tout borélien de $T$ est un
borélien de $T$.

\begin{prop}\label{prop:sat-bor}
  Si $A$ est un borélien de $T$, alors $sat_{\mathcal R}(A)$ est un
  borélien de $T$.
\end{prop}
\preuve On considère les applications $\pi_L,\pi_R:\mathcal R \to X$
restriction au borélien $\mathcal R$ des projections sur le premier et
deuxième facteur de $X\times X$. Puisque $\mathcal R$ est à classes
dénombrables, ces applications sont à fibres dénombrables et
bi-boréliennes en vertu du théorème \ref{thm:kal}.  On remarque alors
que $sat_{\mathcal R}(A)=\pi_L(\pi_R^{-1}(A))$.  \qed

\medskip Etant donné une relation d'équivalence discrète $\mathcal R$
sur un ensemble $T$, on définit l'application $$
\sharp:T\to \N\cup
\{\infty\} $$
qui assigne à tout $x\in T$ le cardinal de sa $\mathcal
R$-classe d'équivalence.

\begin{prop}\label{prop:car-mes}
  Si $\mathcal R$ est une relation d'équivalence borélienne dis\-crète
  sur un espace de Borel standard $T$, alors $\sharp:T\to \N\cup
  \{\infty\}$ est une application borélienne.
\end{prop}
\preuve On considère pour tout $n\in \N$ l'application
$\pi_1:T^{n+1}\to T$ obtenue par projection sur le premier facteur du
produit. C'est une application borélienne mais en général non
bi-borélienne.  On définit le borélien $$
\mathcal
R^n_*=\{(x_0,\dots,x_n)\in T^{n+1}~|~(x_i,x_j)\in \mathcal R, x_i\neq
x_j\} $$
La restriction de $\pi_1$ à $\mathcal R^n_*$ est à fibres
dénombrables car $\pi_1^{-1}(x)\cap \mathcal R^n_*$ est en bijection
avec l'ensemble des parties finies $\mathcal R(x)$ ayant exactement
$n+1$ éléments.  Par conséquent $\pi_1(\mathcal R^n_*)$, i.e.
l'ensemble des $x\in T$ dont la classe d'équivalence $\mathcal R(x)$ a
au moins $n+1$ éléments, est un borélien de $T$ d'après le théorème
\ref{thm:kal}. On en déduit que $$
\sharp^{-1}(n)=\pi_1( \mathcal
R^{n-1}_*)-\pi_1( \mathcal R^n_*) $$
est un borélien de $T$ pour tout
$n\in \N$.  \qed

\medskip Soit $\mathcal R$ une relation d'équivalence borélienne sur
$T$. L'espace quotient $T/\mathcal R$ est muni de la structure
mesurable représentée par la $\sigma$-algèbre des boréliens $\mathcal
R$-saturés de $T$. Une relation d'équivalence borélienne discrète
$\mathcal R$ sur un espace de Borel standard est dite {\em de type I}
(cf. \cite{CFW}) si le quotient $T/\mathcal R$ est un espace de Borel
standard. Ce n'est pas toujours le cas: par exemple la relation
d'équivalence borélienne discrète engendrée par une rotation d'angle
irrationnel sur le cercle n'est pas de type I. Une conséquence
élémentaire mais remarquable du théorème \ref{thm:kal} est la
suivante:

\begin{prop}
  Soit $\mathcal R$ une relation borélienne discrète sur un espace de
  Borel standard $T$. Alors les trois conditions sont équivalentes:
  \begin{itemize}
  \item $\mathcal R$ est de type I;
  \item $\mathcal R$ admet un domaine fondamental borélien, i.e. il
    existe un borélien $D$ de $T$ qui rencontre chaque classe de
    $\mathcal R$ en exactement un point.
  \item Il existe une application borélienne $\nu:T\to \N$ qui est
    injective sur chaque classe de $\mathcal R$.
  \end{itemize}
\end{prop}

Un peu moins élémentaire est le fait que toute relation
d'équiva\-lence finie est de type I, utilisé ou cité sans
démonstration ni référence par divers auteurs.  S'agissant d'un
résultat crucial dans notre étude, et n'ayant pas trouvé de
démonstration dans la littérature, nous en donnons ici une basée sur
le théorème de Kallman:

\begin{prop}\label{lem:dom-fon}
  Soit $\mathcal R$ une relation d'équivalence borélienne finie sur un
  espace de Borel standard $T$. Alors $\mathcal R$ admet un domaine
  fondamental borélien.
\end{prop}
\preuve On fixe $n\in \N$. On peut supposer en vertu de la proposition
\ref{prop:car-mes} que toutes les classes de $\mathcal R$ ont
exactement $n$ éléments. On plonge l'espace de Borel standard $T$ dans
la droite réelle, celle-ci munie de son ordre total naturel, et on
définit un domaine fondamental $D$ de $\mathcal R$ en fixant le plus
petit point de chaque classe. Nous allons montrer que $D$ est un
borélien de $T$.

Il est bien connu (voir par exemple \cite{Rud}) que l'application
$o:T^{n+1}\to T$ définie par $$
o(x_0,\dots,x_n)=\min\{x_0,\dots,x_n\}
$$
est une application borélienne. La restriction de $o$ au borélien
$\mathcal R^n_*$, défini dans la preuve de la proposition
\ref{prop:car-mes}, est à fibres dénombrables. L'ensemble
$D=o(\mathcal R^n)$ est alors un borélien de $T$ en vertu du théorème
\ref{thm:kal}.  \qed

\subsection{BT-espaces}\label{sec:BT-esp}
On appellera {\em espace mesurable-topologique ou MT-espace} tout
ensemble $X$ muni simultanément d'une structure mesurable ou
$\sigma$-algèbre $\mathcal B$ et d'une topologie $\tau$.  Le
composantes connexes de $X$ par rapport à la topologie $\tau$ seront
appelées {\em feuilles}. Ce sont les classes d'équivalence d'une
relation sur $X$ notée $\mathcal R_X$ que nous appellerons {\em
  relation feuille}. Pour toute partie $A$ de $X$ on notera $sat(A)$
le saturé de $A$ par la relation $\mathcal R_X$, i.e. la réunion des
feuilles de $X$ qui rencontrent $A$. Une application $f:X\to Y$ entre
deux MT-espaces sera dite MT si elle est en même temps continue et
mesurable.

\begin{defn}\label{def:MT-esp}
  Un MT-espace $(X,\mathcal B,\tau)$ sera dit {\em standard} ou {\em
    BT-espace} s'il vérifie les conditions suivantes:
  \begin{enumerate}
  \item L'espace mesurable sous-jacent $(X,\mathcal B)$ est un espace
    de Borel standard;
  \item L'espace topologique sous-jacent $(X,\tau)$ est séparé et à
    compo\-santes connexes (feuilles) séparables;

  \item La relation feuille $\mathcal R_X$ est un borélien du produit
    $X\times X$.
  \end{enumerate}
\end{defn}

Soit $X$ un MT-espace. Un {\em MT-sous-espace} de $X$ est un ensemble
$Y$ de $X$ muni de la structure mesurable et la topologie induites. Si
$X$ est un BT-espace et $Y$ est un borélien de $X$, alors la structure
mesurable induite sur $Y$ est standard.  En revanche, il n'est pas sûr
que la relation feuille $\mathcal R_Y$ soit un borélien de $Y\times
Y$. C'est le cas par exemple si $Y$ est un boréliens saturé. En effet
dans ce cas $sat_Y(A)=sat_X(A)$ pour tout $A\subset Y$, et la
condition $3$ de la définition est alors vérifiée par $Y$. Plus
généralement, si $Y$ est un borélien de $X$ tel que la trace de $Y$
sur chaque feuille de $X$ est connexe, alors $Y$ est un MT-sous-espace
standard de $X$.

\label{def:dec}
On appelle {\em découpage} d'un BT-espace une partition dénombrable
$X_i$ de $X$ par des boréliens saturés. Comme remarqué ci-dessus, si
$X$ est un BT-espace, alors les éléments du découpage sont des
BT-espaces car saturés.

Une {\em transversale} d'un BT-espace $X$ est par définition un
borélien $T\subset X$ qui rencontre chaque feuille de $X$ le long d'un
ensemble fermé discret. Autrement dit, il s'agit d'un borélien discret
fermé de $X$. Une transversale est dite {\em totale} si elle rencontre
toutes les feuilles de $X$.  Une transversale totale de $X$ est un
{\em domaine fondamental} si elle rencontre chaque feuille de $X$ en
exactement un point.

La trace des feuilles de $X$ sur une transversale borélienne $T$
déter\-mine une relation d'équivalence sur $T$ que nous noterons
$\mathcal R_X^T$. C'est une relation d'équivalence borélienne sur $T$
au sens de \S\ref{sec:rel-mes} car $\mathcal R_X^T=\mathcal R_X\cap
(T\times T)$. En plus puisque les feuilles de $X$ sont séparables les
classes d'équivalence de $\mathcal R_X^T$ sont toutes dénombrables.

\subsubsection{Mesures quasi-invariantes.}\label{sec:mes}
Soit $X$ un BT-espace et soit $\mathcal T(X)$ l'ensemble de toutes les
transversales boréliennes de $X$.

Une {\em mesure transverse} sur $X$ est par définition une application
$\sigma$-additive $\mu:\mathcal T(X)\to [0,+\infty]$. On supposera
pour simplifier que $\mu$ est sans atome, i.e. que tout point est de
mesure nulle.  Une transversale $T$ est dite {\em $\mu$-négligeable}
si $\mu(T)=0$.

Plus généralement un borélien $A\subset X$ sera dit $\mu$-négligeable
si toute transversale borélienne de $X$ contenue dans $A$ est
$\mu$-négligeable. Une mesure transverse $\mu$ est dite {\em
  quasi-invariante} si le saturé de tout borélien $\mu$-négligeable
est un borélien $\mu$-négligeable. On définit de façon analogue la
notion de mesure quasi-invariante pour une relation borélienne
discrète $\mathcal R$ définie sur un espace de Borel standard $T$.

Une mesure transverse quasi-invariante $\mu$ sur un BT-espace $X$
détermine par restriction une mesure sur toute transversale borélienne
$T\subset X$ que nous noterons encore $\mu$, qui sera quasi-invariante
pour la relation d'équivalence discrète $\mathcal R_X^T$ induite par
$\mathcal R_X$ sur $T$. Dans ce cas le triplet $(T,\mathcal
R_X^T,\mu)$ constitue une {\em relation d'équivalence discrète
  mesurée} au sens de \cite{CFW} ou \cite{FM}.

\subsubsection{Quelques exemples.}
Un BT-espace est un espace de Borel standard $X$ muni d'une topologie
fine compatible en un sens faible avec la structure borélienne.
Celle-ci n'est pourtant pas engendrée en général par la topologie de
$X$, et c'est pour cette raison que l'on est obligé de tenir compte
des deux structures simultanément.

Des exemples triviaux de BT-espaces sont les espaces de Borel standard
munis de la topologie discrète, dits des BT-espaces {\em discrets}.
Dans ce cas les feuilles de $X$ sont réduites à des points de sorte
que la relation feuille $\mathcal R_X$ est la diagonale de $X\times
X$, qui est un borélien produit. Ce type de BT-espaces est caractérisé
par sa structure mesurable. Remarquons en passant qu'il n'y a donc que
deux BT-espaces discrets infinis: l'intervalle $[0,1]$ et l'ensemble
$\N$.

A l'opposé des BT-espaces discrets on trouve les espaces Polonais
connexes. La topologie d'un espace Polonais engendre par définition
une structure borélienne standard. Un espace Polonais est séparé et
séparable et le fait d'être connexe implique que la relation feuille
est triviale et donc borélienne. Par contre nous ignorons sous quelles
conditions un espace Polonais ayant une quantité non dénombrable de
composantes connexes est un BT-espace.

\subsubsection{Prismes}
Un {\em prisme} est un BT-espace $B\times T$ produit d'un espace
Polonais connexe $B$ et d'un espace de Borel standard $T$. L'espace
$B$ est appelé la {\em base} du prisme et $T$ sa {\em verticale}. Nous
appelons {\em plaques} de $B\times T$ les ensembles de la forme
$B\times \{t\}$ avec $t\in T$. On remarquera que la base $B$ étant
connexe, les plaques sont précisément les feuilles du prisme.

Dans ce travail les prismes que nous utiliserons seront
essentiellement ceux dont la base est soit le disque fermé unité
$\mathbb D$ de dimension $2$, soit son bord le cercle $\mathbb S^1$,
soit son intérieur la boule $\mathbb B$. Pour la notion de
triangulation borélienne décrite plus bas nous serons amenés à
remplacer le disque $\mathbb D$ par le simplexe standard $\Delta$ de
dimension $2$. Nous appelons le prisme correspondant un prisme de
disques, cercles, boules ou triangles.

\begin{lem}\label{lem:ouv-pri}
  Soit $B$ un espace Polonais compact connexe et localement connexe et
  $T$ un espace de Borel standard. Soit $A$ un borélien de $B\times
  T$. Alors on a:
  \begin{itemize}
  \item[i)] l'intérieur $int(A)$ et l'adhérence $\overline A$ de $A$
    sont des boréliens de $B\times T$;

  \item[ii)] $int(A)$ est en plus un BT-sous-espace de $B\times T$.
  \end{itemize}
\end{lem}
\preuve (i) Il suffit de montrer que $U=int(A)$ est un borélien de
$B\times T$. Pour voir que $\overline A$ est borélien on applique ce
fait au complé\-mentaire de $A$. On se fixe une base dénombrable
$\mathcal V=\{V_i~|~i\in I\}$ de $B$ formée par des ouverts connexes.
Une partie finie $C\subset I$ sera dite connexe si $B_C=\cup_{i\in C}
V_i$ est connexe dans $B$. Pour toute partie connexe $C\subset I$ on
note $T_C$ l'ensemble des $t\in T$ tels que $B_C\times \{t\}\subset
T$. Autrement dit on a: $$
T_C=T-\pi_T(U-B_C\times T) $$
Les fibres du
borélien $U-B_C\times T$ pour la projection $\pi_T:B\times T\to T$
sont fermées dans un ouvert de $B$. En particulier elles sont réunion
dénombrable de compacts car $B$ est un espace métrique compact.
D'après le théorème \ref{thm:kal-for} la projection $\pi_T(U-B_C\times
T)$ est un borélien de $T$.  On en aura donc de même pour $T_C$. Nous
obtenons ainsi une famille dénombrable de prismes ouverts
$U_C=V_C\times T_C\subset U$ dont la réunion est $U$. Par conséquent
$U$ est un borélien de $B\times T$.

\medskip\noindent (ii) Puisque les feuilles de $U$ sont des espaces
métrisables, elles sont séparées et séparables. Pour voir que $U$ est
un BT-espace, il suffit donc de montrer que sa relation feuille
$\mathcal R_U$ est un borélien de $U\times U$.  La relation feuille
$\mathcal R_{U_C}$ de chaque $U_C$ est un borélien de $X\times X$ car
$U_C$ est un prisme. Nous avons par ailleurs $$
\mathcal R_U=\bigcup_C
\mathcal R_{U_C} $$
où $C$ parcourt les parties finies connexes de
$I$. On en déduit que $\mathcal R_U$ est un borélien de $U\times U$.
\qed

\subsection{Laminations boréliennes}
Nous introduisons ici l'analogue dans le cadre mesurable des atlas
feuilletés d'une variété topologique (voir par exemple \cite{HH}).

Soit $X$ un ensemble et $T$ un espace de Borel standard. Un {\em atlas
  feuilleté borélien de dimension $n$} sur $X$ est un ensemble
dénombrable de couples $\mathcal A=\{(U_i,\varphi_i)\}$, dites {\em
  cartes feuilletées}, où les $U_i$ forment un recouvrement ouvert
localement fini de $X$ et les $\varphi_i$ sont des bijections: $$
\varphi_i:U_i\to V_i\times T_i $$
où $V_i$ est un ouvert connexe de
$\R^n$ et $T_i$ un borélien de l'intervalle $[0,1]$.  Les changements
de cartes $\varphi_i\circ\varphi_j^{-1}$ sont supposés des
BT-isomorphismes locaux. Les ensembles de la forme
$\varphi_i^{-1}(\R^n\times \{*\})$ sont appelés les {\em plaques} de
$\mathcal A$. Les restrictions des $\varphi_i\circ\varphi_j^{-1}$ aux
plaques donnent des homéomorphismes locaux $f_{ij}(t)$ appelés {\em
  changements de plaque}. On dira que $\mathcal A$ est un atlas de
{\em classe $C^p$} si ses changements de plaque sont des
homéomorphismes locaux de classe $C^p$. Deux atlas feuilletés
boréliens $\mathcal A$ et $\mathcal B$ de $X$ sont dits {\em
  $C^p$-compatibles} si $\mathcal A\cup\mathcal B$ est un atlas de
classe $C^p$ de $X$. Un atlas est dit {\em orienté} si les changements
de cartes préservent l'orientation canonique de $\R^n$.

\begin{defn}
  Une {\em lamination borélienne de dimension $n$ et classe $C^p$} sur
  $X$ est par définition la famille $\mathcal F$ de tous les atlas
  $C^p$-compati\-bles avec un atlas donné. Une lamination est dite
  orientable si elle admet un atlas orienté.
\end{defn}

Une lamination $\mathcal F$ de dimension $n$ sur un ensemble $X$
détermine naturellement sur celui-ci une structure MT. En effet, soit
$\{(U_i,\varphi_i)\}$ un atlas feuilleté borélien de $\mathcal F$. Les
bijections $\varphi_i$ induisent sur les ensembles $U_i$ une structure
MT donnée par identification de $U_i$ et $\varphi_i(U_i)\subset
\R^n\times T_i$. On notera $(X,\mathcal F)$ le MT-espace ainsi défini.
Sa structure MT est bien définie modulo des atlas compatibles.

Remarquons que contrairement au cas topologique, nos atlas sont par
définition dénombrables, et que la réunion arbitraire d'atlas n'est
donc pas en général un atlas. En particulier, il n'existe pas d'atlas
maximal définissant $\mathcal F$. Pour remédier à cela on pourrait
considérer des atlas non dénombrables. Mais dans ce cas les structures
mesura\-bles définies par des atlas compatibles ne coïncident pas
forcément. Considérons par exemple un atlas dénombrable $\mathcal A$
et l'atlas non dé\-nombrable $\mathcal A'$ dont les cartes sont les
restrictions de $\mathcal A$ aux plaques. La structure mesurable
induite par $\mathcal A$ est standard tandis que celle induite par
$\mathcal A'$ est transversalement discrète non dénombrable et par
conséquent non standard.

\medskip La preuve du résultat suivant est similaire à celle de
\ref{lem:ouv-pri}.

\begin{prop}\label{prop:ouv-lam}
  Soit $(X,\mathcal F)$ une lamination borélienne à feuilles séparées.
  Alors $(X,\mathcal F)$ est un BT-espace. L'intérieur de tout
  borélien de $X$ est une lamination borélienne de $X$.
\end{prop}

La condition sur la séparabilité des feuilles est naturelle. Dans le
cas topologique cette condition est garantie par la séparabilité de
l'espace topologique ambiant, mais dans notre cas l'espace ambiant
n'est même pas topologique. Toutes les laminations considérées sont
supposées à feuilles séparées, c'est à dire des BT-espaces.

\section{Triangulations et hypercompacité}
\subsection{Laminations par surfaces triangulées}
On supposera désormais que toutes les laminations boréliennes sont de
dimension $2$.

\subsubsection{Piles.}
Soit $(X,\mathcal F)$ une lamination borélienne par surfaces. Une {\em
  pile} de $(X,\mathcal F)$ est donnée par un borélien $\Pi\subset X$
muni d'un BT-isomorphisme: $$
\pi:\Omega\times T\to\Pi $$
où
$\Omega\times T$ est un prisme de base une surface connexe à bord
$\Omega$.  On appelle $\Omega$, $T$ et $\pi$ respectivement la {\em
  base}, la {\em verticale} et le {\em paramétrage} de la pile $\Pi$.
Les plaques de $\Pi$ sont les ensembles $\Pi_t=\pi(\Omega\times
\{t\})$. Chaque plaque de $\Pi$ est munie d'un homéomorphisme ou
paramétrage $\pi_t:\Omega\to \Pi_t$ obtenu par restriction de $\pi$.
L'application $\pi_t$ dépend de façon borélienne de $t$ dans un sens
évident.

\medskip Une pile est donc un borélien de $\Pi$ muni d'une structure
de BT-espace produit induite par le paramétrage. Une pile admet
beaucoup de paramétrages induisant la même structure produit, si l'on
considère celle-ci comme étant la simple donnée des verticales et des
horizontales. Dans le but de fixer cette idée, on dira que deux
paramétrages $$
\pi:\Omega\times T\to \Pi\quad,\quad
\pi':\Omega'\times T'\to \Pi $$
d'une même pile $\Pi$ sont {\em
  équivalents} s'il existe un homéomorphisme $f:\Omega'\to \Omega$ et
un isomorphisme borélien $\gamma:T'\to T$ tels que $$
\pi^{-1}\circ
\pi'(x,t)=(f(x),\gamma(t)), $$
c'est-à-dire qu'ils définissent sur le
borélien $\Pi$ la même structure de BT-espace produit. Si $\Omega$ et
$\Omega'$ sont munies d'une structure additionnelle (différentiable,
simpliciale...) alors on dira que les paramé\-trages sont équivalents
relativement à la dite structure si celle-ci est préservée par
l'homéomorphisme $f$.

\subsubsection{Triangulations de $(X,\mathcal F)$.}\label{sec:tri}
On définit ici la notion de triangulation d'une lamination borélienne.
Il s'agit de donner uns sens à l'idée vague de famille "borélienne" de
triangulations.

\medskip Une {\em triangulation} $\K$ d'une lamination borélienne
$(X,\mathcal F)$ est donnée par un recouvrement de $X$ par une famille
dénom\-brable de piles $\{\Sigma_i\}_{i\in I}$ de base le $2$-simplexe
standard $\Delta$. On appellera {\em triangles} les plaques des piles
$\Sigma_i$. On demandera que si $\sigma$ et $\tau$ sont deux triangles
avec $\tau\cap \sigma\neq \varnothing$ alors:
\begin{itemize}
\item soit $\sigma\cap\tau$ est un sommet;

\item soit $\sigma\cap\tau$ est une arête et le changement de
  paramétrage $\pi_\sigma^{-1}\circ \pi_\tau$ est une bijection affine
  entre deux arêtes de $\Delta$.
\end{itemize}

Il revient au même de dire que l'intersection de deux piles
$\Sigma_i\cap\Sigma_j$ admet un découpage en piles de sommets ou
d'arêtes (voir \S\ref{sec:BT-esp} pour la définition de découpage). Il
est facile à voir que la deuxième condition est indépendante des
paramétrages des piles $\Sigma_i$, pourvu qu'ils soient choisis dans
la même classe d'équivalence simpliciale. On appelle
$\{\Sigma_i\}_{i\in I}$ une {\em famille génératrice} de la
triangulation.

\medskip Une triangulation de $(X,\mathcal F)$ détermine une
triangulation au sens classique de chaque feuille, et en particulier
de l'espace topologique sous-jacent au BT-espace $(X,\mathcal F)$.
Nous notons respectivement $\K^2$, $\K^1$ et $\K^0$ l'ensemble de tous
les triangles, arêtes et sommets de la triangulation. Chaque $\K^i$
peut être identifié à une transversale totale de $X$ formée par les
barycentres des simplexes de $\K^i$. Il hérite donc d'une structure
d'espace de Borel standard. L'ensemble $\K^0\cup\K^1\cup\K^2$ sera
noté $\K$ dans ce papier. Ceci introduit une ambiguïté dans la
notation, car on appelle de la même façon la triangulation (donnée par
une famille génératrice) et le complexe simplicial abstrait
sous-jacent. Cette ambiguïté, très courante dans la littérature, a
comme but d'alléger la notation.

\medskip Si $\mathcal A$ est un sous-ensemble de $\mathcal K$, on note
$|\mathcal A|$ le polytope réunion des simplexes de $\mathcal A$.
L'espace $|\mathcal A|$ est compact si et seulement si $\mathcal A$
est fini.

\begin{lem}\label{lem:pol-bor}
  Soit $A$ un borélien de $\K$. On a alors:
  \begin{itemize}
  \item[i)] $|A|$ est un borélien de $X$. En particulier $|\K^0|$ et
    $|\K^1|$ sont des boréliens de $X$;

  \item[ii)] $sat_X(A)=sat_X(|A|)$ est un borélien de $X$
  \end{itemize}
\end{lem}
\preuve (i) On peut supposer que $A$ est contenu dans une pile
génératrice $\varpi:\Delta\times T\to \Sigma$ de $\K$. Pour tout
simplexe $\tau \subset \Delta$, on note $A_\tau =\varpi(\tau\times
T)\cap A$, qui est borélien car $\varpi$ est un BT-isomorphisme. On
considère la projection $\pi_T$ de $\Sigma$ sur sa verticale $T$. La
restriction de $\pi_T$ à $A$ est à fibres dénombrables, donc
bi-borélienne d'après \ref{thm:kal}. Les ensembles
$T_\tau=\pi_T(A_\tau)$ et $|A_\tau|=\varpi(\tau \times T_\tau)$ sont
donc boréliens. Ceci suffit pour conclure car $|A|=\cup_{\tau}
|A_\tau|$ où $\tau$ parcourt les simplexes de $\Delta$.

\medskip
\noindent (ii) On considère la relation d'équivalence mesurable
discrète $\mathcal R_X^\K$ induite par la relation feuille $\mathcal
R_X$ sur la transversale borélienne $\K$. D'après la proposition
\ref{prop:sat-bor}, l'ensemble $sat_{\mathcal R_X^\K}(A)$ est un
borélien de $\K$. L'ensemble $sat_X(A)=sat_X(|A|)=|sat_{\mathcal
  R_X^\K}(A)|$ est un borélien de $X$ d'après (i).  \qed

\medskip La proposition suivante sera démontrée dans la section
\ref{sec:exi-tri}:

\begin{prop}\label{prop:exi-tri}
  Toute lamination borélienne de classe $C^1$ admet une triangulation.
\end{prop}

Nous supposerons désormais que toute lamination borélienne est munie
d'une triangulation, fixée une fois pour toutes.

\subsubsection{Piles simpliciales.}
Soit $(X,\mathcal F)$ une lamination borélienne. Un borélien $B\subset
X$ est dit {\em simplicial} si toutes ses feuilles sont réunion de
triangles. On appelle {\em volume} d'une feuille de $B$ le nombre de
triangles qu'elle contient. On dira que $B$ est de {\em type fini} si
toutes ses feuilles sont de volume fini.  Etant donné un borélien
simplicial $B$ on note $$
vol:B\to \N\cup\{\infty\} $$
la fonction qui
assigne à chaque $x\in B$ le volume de la feuille de $B$ passant par
$x$. C'est une fonction borélienne d'après la proposition
\ref{prop:car-mes}.

\begin{rem}\label{rem:vol-1-pil}
  Tout borélien $B$ de volume $1$ est une pile. En effet, si
  $\{\Sigma_i\}$ est une famille génératrice de la triangulation,
  alors $B_i=B\cap \Sigma_i$ est un découpage de $B$ en piles. Le
  borélien $B$ est donc une réunion disjointe de piles, et donc une
  pile.
\end{rem}

\medskip Une pile $\Pi$ de $(X,\mathcal F)$ sera dite {\em
  simpliciale} si son paramétrage est simplicial. En particulier $\Pi$
est un borélien simplicial. Les plaques d'une pile étant par
définition des sous-surfaces des feuilles de $(X,\mathcal F)$, si
$\Pi$ est de type fini, alors ses plaques sont réunion d'un nombre
fini de triangles.

\subsubsection{Applications semi-simples.}\label{def:sem-sim}
Soit $B$ un borélien simplicial de $(X,\mathcal F)$ à feuilles
compactes et $\mathcal E$ un ensemble. Une application $f:B\to
\mathcal E$ est dite {\em semi-simple} s'il existe un découpage de $B$
en piles simpliciales $\Pi_i$ tel que pour tout $j$ on a un diagramme
commutatif: $$
\xymatrix{
  \Pi_j \ar[r]^{\rho} \ar[d]_{p_j} & \mathcal E\\
  \Omega_j \ar[ur]_{f_j} } $$
où $p_j$ est la projection de $\Pi_j$
sur sa base $\Omega_j$ et $f_j$ une application de $\Omega_j$ sur $S$.
L'application $f$ est dite {\em simple} s'il existe un $j$ tel que
$B=\Pi_j$.

\begin{lem}\label{lem:sim-sem}
  Soit $\Pi$ une pile simpliciale à base compacte de $(X,\mathcal F)$
  et $K$ un complexe simplicial dénombrable. Alors toute
  BT-application simpliciale $f:\Pi\to K$ est semi-simple.
\end{lem}
\preuve On considère le paramétrage $\pi:\Omega\times T\to \Pi$ de
$\Pi$. Pour tout $t\in T$ on notera $f_t:\Omega\to K$ l'application
simpliciale obtenue par restriction de $f$ à la plaque $\Omega\times
\{t\}$. Etant donnée une application simpliciale $g:\Omega\to K$ on
pose $$
T_g=\{t\in T~|~f_t=g\} $$
Alors $T_g$ est un borélien de $T$.
En effet, $f$ et $g$ étant simpliciales, elles sont caractérisées par
leur valeurs sur les sommets de $\Omega$. Or pour tout sommet $v\in
\Omega$ l'ensemble $$
T_{v,g}=\{t\in T~|~f(v)=g(v)\} $$
est un
borélien car $T_{v,g}\times \{v\}=f^{-1}(g(v))\cap T\times \{v\}$. On
aura alors $$
T_g=\cap_{v\in \Omega} T_{v,g}.  $$

Enfin, puisque $\Omega$ est compacte et $K$ dénombrable, il n'y a
qu'une quantité dénombrable d'applications simpliciales entre eux. La
famille de boréliens $T_g$ détermine alors un découpage $\Pi_g$ de la
pile $\Pi$ vérifiant par construction les conditions requises.  \qed

\medskip Pour généraliser ce résultat à tout borélien simplicial
quelconque de $(X,\mathcal F)$, il faut au préalable savoir découper
un tel borélien en piles.

\subsubsection{Découpage en piles d'un borélien simplicial.}
Le but de ce paragraphe est la preuve du théorème suivant:

\begin{thm}\label{thm:dec-pil}
  Soit $B$ un borélien simplicial de $(X,\mathcal F)$ à feuilles
  compactes. Alors il admet un découpage (cf. page \pageref{def:dec})
  en piles simpliciales $\Pi_i$.
\end{thm}

Comme corollaire de ce théorème et du lemme \ref{lem:sim-sem} nous
avons le résultat suivant:

\begin{cor}\label{cor:sim-sem}
  Soit $B$ un borélien simplicial de $(X,\mathcal F)$ à feuilles
  compactes et $K$ un complexe simplicial dénombrable. Alors toute
  application simpliciale $f:B\to K$ est semi-simple.
\end{cor}

Pour la preuve du théorème, on remarque d'abord qu'on peut se ramener
au cas ou toutes les feuilles de $B$ ont le même volume, car la
fonction $vol$ est borélienne et détermine donc un découpage de $B$ en
boréliens de ce type.

\bigskip\noindent {\bf Un cas particulier}\vspace{1ex}

Nous commençons par le cas particulier d'un borélien qui est réunion
de deux piles simpliciales:

\begin{lem}\label{lem:dec-pil-1}
  Soit $B$ un borélien simplicial de volume fini constant. Soient
  $\Pi'$ et $\Pi''$ deux piles simpliciales telles que toute feuille
  de $B$ est la réunion d'exactement une plaque de $\Pi'$ et une
  plaque de $\Pi''$. Alors $B$ admet un découpage en piles
  simpliciales.
\end{lem}

\preuve La condition sur l'intersection des plaques implique que l'on
peut identifier l'espace de feuilles de $B$ et les verticales des
piles $\Pi'$ et $\Pi''$ à un seul et même espace standard $T$. Soient
$\pi':\Omega'\times T \to \Pi'$ et $\pi'':\Omega''\times T \to \Pi''$
les paramétrages de $\Pi'$ et $\Pi''$ respectivement. On considère la
BT-application simpliciale: $$
\pi'\sqcup \pi'':(\Omega'\sqcup
\Omega'')\times T\to B $$
somme disjointe des applications $\pi$ et
$\pi'$. C'est une application bi-borélienne car borélienne et à fibres
dénombrables (cf. théorème \ref{thm:kal}). Pour tout $t\in T$
l'application $\pi'_t\sqcup \pi''_t:\Omega'\sqcup \Omega''\to B_t$
induit sur $\Omega'\sqcup \Omega''$ une relation d'équivalence notée
$\rho_t$. On montre alors que l'ensemble $T_\rho=\{t\in
T~|~\rho_t=\rho\}$ est un borélien de $T$, par une preuve analogue à
celle du lemme \ref{lem:sim-sem}.

\medskip On note $\Omega_\rho$ le complexe simplicial quotient de
$\Omega'\sqcup \Omega''$ par la relation d'équivalence simpliciale
$\rho$, et soit $$
q_\rho:(\Omega'\sqcup \Omega'') \times T_\rho \to
\Omega_\rho\times T_\rho $$
la BT-application quotient, qui est
bi-borélienne car à fibres dénom\-brables.  Nous obtenons ainsi des
diagrammes commutatifs d'applica\-tions bi-boréliennes $$
\xymatrix{
  (\Omega'\sqcup \Omega'')\times T_\rho \ar[rr]^{\qquad
    \pi'\sqcup\pi''}
  \ar[d]_{q_\rho} & & B\\
  \Omega_\rho\times T_\rho \ar[rru]^{\pi_\rho} } $$
L'application
$\pi_\rho$ est injective par construction. C'est donc un
BT-isomorphisme. Le borélien saturé $B_\rho\subset B$, image de
$\pi_\rho$, est en particulier une pile simpliciale. Nous obtenons
ainsi un découpage de $B$ en piles simpliciales $B_\rho$, où $\rho$
parcourt les relations d'équivalence simpliciales sur le complexe
$\Omega'\sqcup \Omega''$.  \qed

\bigskip\noindent {\bf Le cas général}\vspace{1ex}

On montre ici comment réduire le cas général au cas présenté dans le
paragraphe précédent. Pour ce faire nous montrerons que tout borélien
simplicial peut être découpé en boréliens qui vérifient les conditions
du lemme \ref{lem:dec-pil-1}, i.e. qui sont réunion de deux piles
simpliciales.

\begin{lem}\label{lem:dec-pil-2}
  Tout borélien simplicial $B$ admet un découpage en boré\-liens
  simpliciaux vérifiant les conditions du lemme \ref{lem:dec-pil-1}.
\end{lem}
\preuve On peut toujours supposer que les feuilles de $B$ sont de
volume constant. Nous allons procéder par récurrence sur ce volume.
Nous avons déjà vu (cf. remarque \ref{rem:vol-1-pil}) que tout
borélien de volume $1$ est une pile.

Soit $n$ un entier $\geq 1$. Supposons que le lemme
\ref{lem:dec-pil-1} est vérifié par tout borélien simplicial de volume
$\leq n$. Supposons $B$ de volume constant $n+1$. On note $\hat
B=\K^2\cap B$ le borélien des triangles de $B$. Il existe d'après le
lemme \ref{lem:dom-fon} un borélien de $\mathcal J\subset \hat B$ qui
rencontre chaque feuille de $\Pi$ en un seul point.

Le borélien $|\mathcal J|$ est de volume $1$ et donc une pile
simpliciale. Le borélien $|\hat B-\mathcal J|$ est pour sa part de
volume $n$ et admet par hypothèse de récurrence un découpage en piles
simpliciales $\hat B_i$. On note $B_i$ la réunion des feuilles de $B$
rencontrant $\hat B_i$. Chaque feuille de $B_i$ est le réunion d'une
plaque de la pile $\hat B_i$ et d'un triangle de $|\mathcal J|$ et il
vérifie donc le lemme \ref{lem:dec-pil-1}. Enfin on remarque que
puisque $B$ est de volume $\geq 2$, toute feuille de $B$ rencontre
$|\hat B-\mathcal J|$. La réunion des $B_i$ est en particulier égale à
$B$. Ceci complète la récurrence.  \qed

\subsection{Hyperfinitude et hypercompacité}

\subsubsection{Relations d'équivalence et laminations
hyperfinies.}\label{sec:hyp}
Une relation d'équivalence borélienne discrète $\mathcal R$ sur un
espace de Borel standard $T$ est dite {\em hyperfinie} si elle est
limite inductive de relations d'équi\-valence boréliennes finies. Plus
précisément, il existe une suite croissante $\mathcal R_n$ de
relations d'équivalence boréliennes finies sur $T$ dont la réunion est
$\mathcal R$. On notera $\mathcal R=\varinjlim \mathcal R_n$. On
remarquera que si $\mathcal R$ admet un domaine fondamental borélien,
alors elle est hyperfinie.

Rappelons que si $T$ est une transversale borélienne de $(X,\mathcal
F)$, la trace des feuilles de $X$ sur $T$ détermine une relation
d'équivalence borélienne discrète que nous notons $\mathcal F_T$. La
définition suivante est due à Bowen.

\begin{defn}[\cite{Bow}]
  On dira que $(X,\mathcal F)$ est une lamination {\em hyperfinie} si
  et seulement si la relation d'équivalence $\mathcal F_T$ est
  hyperfinie quelle que soit la transversale $T$ choisie.
\end{defn}

\paragraph{Exemples de feuilletages hyperfinis.}
Soit $V^{n+d}$ et $V^n$ deux variétés de dimension $n+d$ et $n$
respectivement.  Les fibres d'une submersion $p:V^{n+d}\to V^n$
définissent un feuilletage $\mathcal F$ de dimension $d$ sur
$V^{n+d}$. Ce feuilletage est hyperfini. En effet, soit $T$ une
transversale borélienne de $\mathcal F$. La restriction de la
submersion $p$ à $T$ est une application borélienne à fibres
dénombrables qui, d'après le théorème de Kallman, admet une section
borélienne. L'image de cette section est un domaine fondamental pour
la relation $\mathcal F_T$, qui est donc hyperfinie.

Le premiers exemples de feuilletages hyperfinis sans domaine
fondamental sont les feuilletages linéaires du tore $\mathbb T^3$,
définis par l'équa\-tion différentielle $dz=\alpha dx + \beta dy$,
avec $\alpha,\beta\in \R$. Ils sont donnés par une action
différentiable de $\R^2$, et sont donc hyperfinis d'après \cite{JKL}.
Si $\alpha$ et $\beta$ sont rationnellement indépendants, alors
l'action de $\R^2$ est libre et les feuilles du feuilletage sont des
plans.

\paragraph{Un exemple de feuilletage {\em non} hyperfini.}\label{ex:non-hyp}
Soit $\Sigma_2$ la surface de genre $2$. On considère une
représentation fidèle du groupe fondamental $\pi_1(\Sigma_2)$ de
$\Sigma_2$ dans le groupe $SO(3)$. La suspension de l'action libre
correspondante de $\pi_1(\Sigma_2)$ sur $SO(3)$ est une variété
compacte $V$ de dimension $5$ munie d'un feuilletage par plans
$\mathcal F$.  L'action de $\pi_1(\Sigma_2)$ sur $SO(3)$ préserve la
mesure de Haar (finie) de celui-ci et engendre par conséquent une
relation d'équivalence discrète qui n'est pas hyperfinie car le groupe
$\pi_1(\Sigma_2)$ est non moyennable (cf.  \cite{Zim} et \cite{CFW}).
Une autre façon de voir que la relation est non hyperfinie consiste à
prouver que son {\em coût} est égale à $2$ (cf.  \cite{Gab}). Le
feuilletage $\mathcal F$ est donc {\em non hyperfini}.

\subsubsection{Laminations hypercompactes.}
Rappelons qu'on a fixé une triangulation de $(X,\mathcal F)$, et que
l'ensemble des ses simplexes $\K$ est muni d'une structure de Borel
standard donnée par identification de chaque simplexe avec son
barycentre dans $X$. L'espace $\K$ est du coup identifié à une
transversale borélienne de $(X,\mathcal F)$. On note $\hat{\mathcal
  F}$ la relation d'équivalence induite sur l'espace des triangles
$\K^2$ par les feuilles de $X$. La relation $\hat{\mathcal F}$ n'est
autre que la relation $\mathcal R^{\K^2}_X$ introduite dans
\S\ref{sec:BT-esp}. La raison de ce changement de notation sont
claires. Si $(X,\mathcal F)$ est hyperfinie, alors $\hat{\mathcal F}$
sera hyperfinie. La réciproque sera montrée ci-dessous.

\begin{defn}
  Une lamination $(X,\mathcal F)$ sera dite {\em hypercompacte} s'il
  existe une suite $B_n$ de boréliens simpliciaux pour $\K$ de type
  fini telle que
  \begin{enumerate}
  \item $B_n\subset B_{n+1}$ pour tout $n\in \N$;
  \item $\cup_n B_n =X$.
  \end{enumerate}
  La suite $B_n$ sera appelé une {\em filtration compacte simpliciale}
  de $(X,\mathcal F)$.
\end{defn}

\begin{thm}\label{thm:hyp-hyp}
  Pour toute lamination borélienne $(X,\mathcal F)$ les deux
  conditions suivantes sont équiva\-lentes:
  \begin{enumerate}
  \item[a)] $(X,\mathcal F)$ est hypercompacte;
  \item[b)] $(X,\mathcal F)$ est hyperfinie.
  \end{enumerate}
\end{thm}

\subsubsection{Preuve de (a)$\Rightarrow$(b).}
Soit $B_n$ une filtration compacte de $(X,\mathcal F)$. Soit $T$ une
transversale borélienne et soit $\mathcal F_T$ la relation induite par
$\mathcal F$ sur $T$. Puisque les $B_n$ sont à feuilles compactes, la
relation $\mathcal R_n$ induite par $\mathcal F$ sur la transversale
borélienne $\hat B_n=T\cap B_n$ est une sous-relation borélienne finie
de $\mathcal F_T$. Puisque $B_n\subset B_{n+1}$ on a $\mathcal
R_n\subset \mathcal R_{n+1}$. Enfin du fait que $\cup_n B_n =X$ on en
déduit que $\cup_n \mathcal R_n=\mathcal F_T$. La relation $\mathcal
F_T$ est donc hyperfinie. Ceci complète la preuve de la première
implication.

\subsubsection{Borélien associé à une relation d'équivalence borélienne finie.}
Avant d'aborder la preu\-ve de l'implication réciproque, on aura
besoin de la construction suivante. On se fixe une sous-relation
d'équivalence borélienne finie $\mathcal R$ de $\hat{\mathcal F}$.
Chaque classe d'équivalence $\mathcal R(\sigma)$ de $\mathcal R$
détermine un domaine simplicial compact $|\mathcal R(\sigma)|$ sur la
feuille $L_\sigma$ obtenu par réunion des triangles de $\mathcal
R(\sigma)$.

\medskip On appelle {\em $0$-squelette} de $\mathcal R$ l'ensemble
$sq^0( \mathcal R)\subset \K^0$ formé par les sommets de $\K$ qui sont
dans la frontière d'un des domaines $|\mathcal R(\sigma)|$. On appelle
{\em $1$-squelette} de $\mathcal R$ l'ensemble $sq^1( \mathcal
R)\subset \K^1$ formé par les arêtes de $\K$ dont les sommets sont
dans $sq^0( \mathcal R)$.

\begin{lem}
  Pour $i=0,1$ le $i$-squelette $sq^i(\mathcal R)$ d'une relation
  d'équi\-valence borélienne finie $\mathcal R\subset \hat{\mathcal
    F}$ est un borélien de $\K^i$.
\end{lem}
\preuve On plonge d'abord $\K^2$ dans l'intervalle $[0,1]$ et on
plonge ensuite $\K^1$ dans $\K^2\times\K^2$ en assignant à toute arête
la paire de triangles dont elle est la face commune, ceux-ci étant
ordonnés de façon croissante par rapport à l'ordre total de $[0,1]$.
Une arête est dans le $1$-squelette de $\mathcal R$ si et seulement si
les deux triangles dont elle est face appartiennent à des $\mathcal
R$-classes d'équivalence différentes. On a alors $sq^1( \mathcal
R)=\K^1 - \mathcal R$.

De façon analogue on plonge $\K^1$ dans $\K^0\times \K^0$ en assignant
à toute arête la paire ordonnée de ses sommets. L'ensemble
$sq^1(\mathcal R)$ est donc un borélien de $\K^0\times \K^0$. Soit
$\pi_1,\pi_2:\K^0\times \K^0\to \K^0$ les projections sur le premier
et deuxième facteur. La restriction de chaque $\pi_i$ au borélien
$sq^1(\mathcal R)$ est à fibres dénombrables, d'où l'on déduit que les
$\pi_i( sq^1(\mathcal R))$ sont des boréliens $\K^0$. On remarque
enfin que $$
sq^0(\mathcal R)=\pi_1( sq^1(\mathcal R))\cup \pi_2(
sq^1(\mathcal R)), $$
ce qui complète la preuve du lemme.  \qed

\medskip Etant donné une relation d'équivalence borélienne finie
$\mathcal R\subset \hat{\mathcal F}$, on dira qu'un triangle
$\sigma\in \K^2$ est {\em intérieur} à $\mathcal R$ s'il ne rencontre
pas $|sq^1(\mathcal R)|$. On note $int(\mathcal R)\subset \K^2$
l'ensemble des triangles de $\K$ intérieurs à $\mathcal R$.

\begin{lem}
  L'ensemble $int(\mathcal R)$ est un borélien de $\K^2$.
\end{lem}
\preuve Il suffit de prouver que, étant donnée une pile géné\-ratrice
$\Sigma\simeq \Delta^2\times T$ de $\K$, l'ensemble des triangles de
$\Sigma$ qui sont intérieurs à $\mathcal R$ est un borélien de $T$.
Soient $v_i$ les trois sommets de $\Delta^2$ et soit $T_i$ le borélien
de $T$ défini par $sq^0(\mathcal R) \cap (\{v_i\}\times T )=
\{v_i\}\times T_i$. Alors $T\cap int(\mathcal R)=T_1\cap T_2\cap T_3$
est un borélien.  \qed

\medskip Le borélien de type fini $$
B(\mathcal R)=|int(\mathcal R)|
$$
est appelé le {\em borélien associé à} $\mathcal R$. Il est clair
que si $\mathcal R_1\subset \mathcal R_2$ alors $B(\mathcal
R_1)\subset B(\mathcal R_2)$.

\subsubsection{Preuve de (b)$\Rightarrow$(a).}
On considère une suite croissante de relations d'équivalence
boréliennes finies $\mathcal R_n\subset \hat{\mathcal F}$ telle que
$\cup_n\mathcal R_n=\hat{\mathcal F}$. On pose $$
B_n=|B(\mathcal
R_n)|\quad,\quad n\in \N $$
où $B(\mathcal R_n)$ est le borélien de
type fini construit dans le paragraphe précédent. Puisque la suite
$\mathcal R_n$ est croissante, il en est de même pour la suite de
boréliens simpliciaux $B_n$. Il reste à montrer qu'il s'agit d'une
suite exhaustive. Pour ce faire on considère un triangle quelconque
$\sigma\in \K^2$. Puisque $\sigma$ est intérieur à la feuille
$L_\sigma$ et la suite de compacts simpliciaux $|\mathcal
R_n(\sigma)|$ ($n\in \N$) est croissante et exhaustive, il existe un
$n\in \N$ tel que $\sigma\subset int(|\mathcal R_n(\sigma)|)$. En
particulier $\sigma\subset B_n$. Ceci complète la preuve de la
deuxième implication et du théorème.

\section{Laminations par plans}
On supposera désormais que $(X,\mathcal F)$ est une lamination par
plans, c'est à dire une lamination dont toutes les feuilles sont
homéomorphes à $\R^2$. Elle est munie d'une triangulation $\K$.

\subsection{Hypercompacité forte}
Une lamination borélienne $(X,\mathcal F)$ est dite {\em fortement
  hypercompacte} si elle admet une filtration compacte simpliciale
$B_n$ formée par des boréliens dont toutes les feuilles sont des
disques. Le but de cette section est la preuve du théorème suivant:

\begin{thm}\label{thm:hyp-hyp-for} 
  Soit $(X,\mathcal F)$ une lamination borélienne par plans. Alors les
  deux conditions sont équivalentes:
  \begin{itemize}
  \item[a)] $(X,\mathcal F)$ est fortement hypercompacte;
  \item[b)] $(X,\mathcal F)$ est hyperfinie.
  \end{itemize}
\end{thm}

Tout comme dans la preuve de \ref{thm:hyp-hyp}, l'enjeu est de montrer
(b) $\Rightarrow$ (a). Nous allons en fait montrer que
l'hypercompacité implique, dans le cas des laminations par plans,
l'hypercompacité forte. L'idée est simple: il s'agit de prendre une
filtration compacte $B_n$ de $(X,\mathcal F)$ pour ensuite "boucher
les trous" des feuilles des boréliens $B_n$ de façon à obtenir une
nouvelle suite de boréliens simpliciaux par disques $\hat B_n$. Nous
clarifions dans le paragraphe suivant le sens précis de l'expression
"boucher les trous".

\subsection{Enveloppes}\label{sec:env}
\subsubsection{Remarques préliminaires}
\begin{rem}\label{rem:env-1}
  Soit $\Omega$ une surface compacte connexe plongée dans $\R^2$. On
  appelle {\em enveloppe de $S$} le disque topologique
  $\epsilon(\Omega)\subset \R^2$ défini par $$
  \epsilon(\Omega)\supset
  \Omega\quad,\quad \partial\epsilon(\Omega)\subset \partial \Omega $$
  C'est le disque complémentaire de la seule composante connexe non
  bornée de $\R^2-\Omega$. Il est obtenu en rajoutant à $\Omega$ ses
  "trous", i.e. les composantes bornés du complémentaire. Le lecteur
  remarquera que {\em si $\Omega_1$ et $\Omega_2$ sont deux surfaces
    disjointes ou emboîtées dans $\R^2$, alors les disques
    $\epsilon(\Omega_1)$ et $\epsilon(\Omega_2)$ sont disjoints ou
    emboîtés.} En effet, les composantes connexes non bornées de
  $\R^2-\Omega_1$ et $\R^2-\Omega_2$ ne peuvent pas être disjointes.
  Si elles sont emboîtées alors on a de même pour $\epsilon(\Omega_1)$
  et $\epsilon(\Omega_2)$. Sinon on aura
  $\epsilon(\Omega_1)\cap\epsilon(\Omega_2)=\varnothing$.
\end{rem}

\begin{rem}\label{rem:env-2}
  Supposons que $\R^2$ est muni d'une triangulation et soit $\Omega_i$
  une famille dénombrable de surfaces compactes simpliciales deux à
  deux disjointes dans $\R^2$. Alors {\em l'ensemble
    $\hat\Omega=\bigcup_i\epsilon(\Omega_i)$ est soit une réunion de
    disques deux à deux disjoints, soit $\R^2$ tout entier}. En effet,
  d'après la remarque précédente les disques simpliciaux
  $\epsilon(\Omega_i)$ sont deux à deux disjoints ou emboîtés. Puisque
  les surfaces $\Omega_i$ sont disjointes, les bords de leurs
  enveloppes ne se rencontrent pas. Or il est clair que toute suite
  croissante de disques simpliciaux emboîtes à bords disjoints est
  soit finie, soit exhaustive. Si la famille de disques
  $\epsilon(\Omega_i)$ contient des suites croissantes non finies,
  alors $\hat \Omega=\R^2$. Dans le cas contraire $\hat \Omega$ est
  une réunion disjointe de disques simpliciaux, à savoir les éléments
  maximaux de la famille $\epsilon(\Omega_i)$.
\end{rem}

\subsubsection{Enveloppe d'un borélien simplicial de $(X,\mathcal F)$.}
Soit $B$ un boré\-lien simplicial dont toutes les feuilles sont des
surfaces compactes. Soit $T$ l'espace de ses feuilles. Pour tout $t\in
T$ on note $B_t$ la feuille de $B$ correspondante. On définit
l'enveloppe de $B$ comme l'ensemble $$
\epsilon(B)=\bigcup_{t\in T}
\epsilon(B_t) $$
où $\epsilon(B_t)$ est l'enveloppe de la surface
compacte $B_t$ dans la feuille de $(X,\mathcal F)$ qui la contient.

\begin{lem}
  L'enveloppe $\epsilon(B)$ d'un borélien simplicial à feuilles
  compactes $B$ est un borélien simplicial.
\end{lem}
\preuve On considère le borélien ouvert $$
C=X-B $$
complémentaire de
$B$ dans $X$. L'adhérence de chaque feuille de $C$ est une surface
éventuellement non compacte. On note $C_\infty$ le borélien saturé de
$C$ formé par la réunion de ses feuilles non bornées, i.e. non
relativement compactes. C'est un borélien en vertu de la proposition
\ref{prop:car-mes}. De plus:
\begin{equation}\label{eq:env}
  \epsilon(B)=X-C_\infty.
\end{equation}
En effet, soit $L$ une feuille de $(X,\mathcal F)$ et soit $\Omega_i$
la famille des feuilles de $B$ contenues dans $L$. Puisque $L$ a un
seul bout, $L-B$ a au plus une composante connexe non bornée $U$,
intersection des composantes connexes non bornées de $L-\Omega_i$. Le
complémentaire de $U$ dans $L$ est alors la réunion des enveloppes des
surfaces $\Omega_i$. Autrement dit, on a: $$
L-C_\infty=L-U=\cup_i
\epsilon(\Omega_i)=L\cap \epsilon(B).  $$
Ceci montre l'identité
(\ref{eq:env}) et donc le lemme.  \qed

\medskip Le lecteur remarquera que le borélien $\epsilon(B)$ n'est pas
forcément à feuilles compactes. Il se peut même que $\epsilon(B)=X$;
imaginons un borélien dont la trace sur chaque feuille est une famille
de couronnes concentriques disjointes. Pour contourner cette
difficulté on introduit la notion de borélien intègre.

On dira que $B$ est un borélien {\em intègre} si pour toute feuille
$\Omega$ de $B$, on a $\epsilon(\Omega)\cap B= \Omega$. Autrement dit,
les enveloppes des feuilles de $B$ sont deux à deux disjointes. Ainsi
il existe une identification naturelle entre les feuilles de $B$ et
les feuilles de $\epsilon(B)$ de sorte que toute feuille de
$\epsilon(B)$ est l'enveloppe de la feuille de $B$ correspondante. En
particulier, $\epsilon(B)$ est un borélien de disques. Nous montrons
d'abord le lemme suivant.

\begin{lem}\label{lem:dec-int}
  Tout borélien simplicial $B$ admet un découpage en boré\-liens
  intègres.
\end{lem}
\preuve L'idée est très simple: le borélien $B$ contient des feuilles
"minimales" dont l'enveloppe ne contient aucune autre feuille de $B$.
On note $B_*$ le saturé de $B$ formé par ces feuilles-là. C'est un
ensemble non vide par le lemme de Zorn.  Le problème est de montrer
qu'il s'agit d'un borélien. Une fois montré ceci on procède par
récurrence en posant $$
B^{i+1}=B^i-B^i_*\quad,\quad B^0=B $$
Il est
clair par construction que les boréliens $B^i_*$ sont tous intègres.
Il est aussi facile à voir que la suite $B^i_*$ est exhaustive. En
effet, si l'enveloppe d'une feuille donnée $B_t$ contient exactement
$k$ autres feuilles de $B$, alors $B_t$ appartient à $B^k_*$.

Pour montrer que $B_*$ est borélien, on considère $\overline C$
l'adhérence dans $X$ du borélien $C=X-B$. C'est le complémentaire dans
$X$ de $int(B)$, qui est un borélien d'après la proposition
\ref{prop:ouv-lam}. Soit $E$ l'ensemble saturé de $\overline C$ formé
par les feuilles compactes de $\overline C$. Il admet un découpage en
piles d'après le théorème \ref{thm:dec-pil}. Soit $E_*$ la réunion des
piles de $E$ dont la base est un disque. Les feuilles de $E_*$ ne sont
que les "trous" des feuilles de $B_*$, et par conséquent $B_*$ est le
saturé dans $B$ de l'ensemble $$
\partial E_* = E_*\cap B.  $$
Or
$\partial E_*$ est un borélien simplicial et le saturé de tout
borélien simplicial est borélien (voir le lemme \ref{lem:pol-bor}). Le
lemme est démontré.  \qed

\medskip Nous avons besoin d'un dernier lemme technique avant de nous
attaquer à la preuve du théorème \ref{thm:hyp-hyp-for}. On se fixe un
borélien simplicial à feuilles compactes $B$ et on reprend la notation
introduite au début de ce paragraphe. On pose $$
T^{(q)}=\{t\in
T~|~vol(\epsilon(B_t))\leq q\}\quad,\quad \epsilon_q(B)=\bigcup_{t\in
  T^{(q)}} \epsilon(B_t).  $$

\begin{lem}\label{lem:env-bor}
  Pour tout $q\in \N$ l'ensemble $\epsilon_q(B)$ est un borélien
  simplicial de $(X,\mathcal F)$ dont toutes les feuilles sont des
  disques.
\end{lem}
\preuve En vertu du lemme \ref{lem:dec-int}, on peut supposer que $B$
est un borélien intègre. Dans ce cas $\epsilon(B)$ est à feuilles
compactes et l'ensemble $\epsilon_q(B)$ est le saturé de $\epsilon(B)$
formé par ses feuilles de volume $q$. C'est un borélien d'après la
proposition \ref{prop:car-mes}.  \qed

\subsection{Preuve du théorème \ref{thm:hyp-hyp-for}}
On se fixe comme annoncé une suite croissante exhaustive $B_n$ de
boréliens simpliciaux de $(X,\mathcal F)$. On définit maintenant pour
tout $q\in N$ le borélien: $$
\hat B_q =\bigcup_{n \in \N}
\epsilon_q(B_n) $$
On veut montrer que la suite $\hat B_q$ est une
filtration compacte forte pour la lamination $(X,\mathcal F)$. Il est
clair que $\hat B_q\subset \hat B_{q+1}$ pour tout $q\in \N$, et il ne
nous reste donc qu'à vérifier:
\begin{enumerate}
\item[i)]$\hat B_q$ est un borélien. c'est le lemme \ref{lem:env-bor}.

\item[ii)] Toute feuille $\Omega$ de $\hat B_q$ est un disque de
  volume $\leq q$; en effet $\Omega$ est réunion d'une suite de
  feuilles $\Omega_{q,n}$ de $\epsilon_q(B_n)$, qui sont, d'après la
  remarque \ref{rem:env-1}, des disques disjoints ou emboîtés de
  volume borné par $q$.

\item[iii)] La suite $\hat B_q$ est exhaustive car $\bigcup_{q\in \N}
  \hat B_q\supset \bigcup_n B_n =X$.
\end{enumerate}
Le théorème est démontré.

\section{Revêtements}
\subsection{Définition}
Soit $S$ une surface compacte. Un {\em BT-revêtement} de $S$ est une
paire $(X,\rho)$ donnée par un BT-espace $X$ et par une
BT-application: $$
\rho:X\to S $$
qui est un revêtement de $S$ en
restriction à chaque feuille de $X$. Dans ce cas $X$ est naturellement
muni d'une lamination borélienne $\mathcal F$ définie par un atlas
feuilleté obtenu par image réciproque d'un atlas fini de $S$. On
définit de façon analogue la notion de BT-homéomor\-phisme local.

Le lecteur remarquera qu'un BT-revêtement n'est rien d'autre qu'un
revêtement au sens classique défini sur l'espace topologique
sous-jacent à $X$, avec la particularité d'être globalement borélien.

\subsubsection{Suspensions.}
Des nombreux exemples de BT-revêtements sont construits par le
pro\-cédé de {\em suspension}, que nous décrivons ici. Soit $\Gamma$
le groupe fondamental de $S$. Supposons qu'il agit de façon borélienne
sur un espace de Borel standard $T$. Nous considérons l'action
diagonale de $\Gamma $ sur le prisme $\tilde X=\tilde S\times T$, où
$\tilde S$ est le revêtement universel de $S$. L'espace quotient de
cette action $X=\tilde X/\Gamma$ est naturellement muni d'une
structure MT dont on vérifie aisément qu'elle est standard. La
projection sur le premier facteur $p_1:\tilde X\to \tilde S$ passe au
quotient en une application de revêtement $\rho:X\to S$, comme
l'indique le diagramme commutatif suivant: $$
\xymatrix{
  \tilde X \ar[r]^{\hat q} \ar[d]_{p_1} & X \ar[d]^{\rho}\\
  \tilde S \ar[r]^q & S } $$
Le BT-revêtement $(X,\rho)$ est appelé la
{\em suspension} de l'action.

Réciproquement, si $(X,\rho)$ est un BT-revêtement de $S$, nous
pouvons présenter $(X,\rho)$ comme une suspension. On se fixe un point
$*\in S$ et on note $T=\rho^{-1}(*)$ la transversale de $X$ préimage
de $*$ par l'application de revêtement $\rho$. Le groupe $\Gamma$ agit
de façon borélienne sur $T$, et la suspension de cette action nous
donne précisément le BT-revêtement $(X,\rho)$ de départ.

\bigskip
\subsection{Laminations hyperfinies et revêtements}
Le but de ce paragraphe est la preuve du résultat suivant:

\begin{thm}\label{thm:rev}
  Soit $(X,\mathcal F)$ une lamination borélienne par plans. Alors les
  trois conditions suivantes sont équivalentes:
  \begin{enumerate}
  \item $(X,\mathcal F)$ est hyperfinie;

  \item $(X,\mathcal F)$ est un BT-revêtement de toute surface de
    genre $\geq 1$;

  \item $(X,\mathcal F)$ est un BT-revêtement de $\mathbb T^2$;

  \item $(X,\mathcal F)$ est définie par une BT-action libre de
    $\R^2$.
  \end{enumerate}
\end{thm}

Les implications (2)$\Rightarrow$(3)$\Rightarrow$(4) sont immédiates.
Le tore $\mathbb T^2$ étant une surface de genre $1$, la condition (3)
est un cas particulier de (2). Si $(X,\mathcal F)$ est un
BT-revêtement de $\mathbb T^2$, alors on peut relever l'action
naturelle de $\R^2$ sur le tore en une BT-action localement libre sur
$(X,\mathcal F)$. Les feuilles de $(X,\mathcal F)$ étant par hypothèse
des plans, cette action n'a pas de point fixe.

\medskip Le lemme suivant, c'est à dire l'implication
(4)$\Rightarrow$(1), est un corollaire d'un résultat de Jackson,
Kechris et Louveau \cite{JKL}, selon lequel une action borélienne d'un
groupe discret moyennable à croissance polynomiale engendre une
relation d'équivalence hyperfinie:

\begin{lem}
  Toute BT-action de $\R^2$ est hyperfinie.
\end{lem}
\preuve Soit $S$ une transversale borélienne de $(X,\mathcal F)$ et on
considère $\hat S$ le saturé de $S$ par l'action borélienne de $\Z^2$
induite sur $X$ par celle de $\R^2$. Le borélien $\hat S$ n'est pas
une transversale car en générale non fermé ou non discret. La relation
d'équivalence borélienne induite par $\mathcal F$ sur $\hat S$ est
néanmoins à classes dénombrables et hyperfinie d'après \cite{JKL}. La
restriction de cette relation au sous-borélien $S$ sera également
hyperfinie.  \qed

\medskip L'implication qui nécessite le plus de travail est
(1)$\Rightarrow$(2), que nous énonçons explicitement dans la
proposition suivante:

\begin{prop}\label{prop:lam-pla-rev}
  Toute lamination par plans hyperfinie est un BT-revêtement de toute
  surface de genre non nul.
\end{prop}

L'idée de la preuve est la suivante. Si $(X,\mathcal F)$ est une
lamination par plans hyperfinie, alors d'après le théorème
\ref{thm:hyp-hyp-for}, elle est fortement hypercompacte, ce qui
signifie qu'elle admet une filtration compacte simplicial $B_n$ formée
par des boréliens dont toutes les feuilles sont des disques. Fixons
une surface de genre non nul $S$. Nous allons construire une
BT-application de revêtement $\rho:X\to S$ comme limite d'une suite de
BT-homéomorphismes locaux $\rho_n:B_n\to S$ définis sur la filtration
$B_n$.

\subsubsection{Découpages relatifs.}
Soit $\Pi$ et $\hat \Pi$ deux piles simpliciales de $(X,\mathcal F)$.
On suppose $\Pi\subset \hat \Pi$.

On dira que $\Pi$ est une {\em sous-pile pleine de} $\hat \Pi$ si
toute plaque de $\hat \Pi$ rencontre $\Pi$ et s'il existe un
plongement simplicial $f:\Omega\to \hat \Omega$ de la base de $\Pi$
dans la base de $\hat \Pi$ et un diagramme commutatif $$
\xymatrix{
  \Pi \ar[r]^i \ar[d]_p & \hat \Pi \ar[d]^{\hat p}\\
  \Omega \ar[r]^f & \hat \Omega } $$
où $p$ et $\hat p$ sont les
projections des piles sur leurs bases respectives.  Les trois
propriétés suivantes découlent immédiatement de la définition de
sous-pile pleine:
\begin{enumerate}\label{pil-ple}
\item[a)] chaque plaque de $\hat \Pi$ contient une et une seule plaque
  de $\Pi$.  En particulier, les verticales de $\Pi$ et $\hat \Pi$
  sont les mêmes;

\item[b)] l'application $\hat p \circ i$ est une application simple au
  sens de la définition \ref{def:sem-sim} (pag.
  \pageref{def:sem-sim}).

\item[c)] Les verticales de $\Pi$ et $\hat \Pi$ coïncident, de sorte
  que tout découpage $\Pi_i$ de $\Pi$ détermine un découpage $\hat
  \Pi_i$ de $\hat \Pi$.
\end{enumerate}

On généralise la définition ci-dessus à des paires de boréliens à
feuilles compactes. On dira qu'un borélien simplicial $B\subset \hat
\Pi$ est un {\em sous-borélien plein} de la pile $\hat \Pi$ s'il admet
un découpage en un nombre fini de sous-piles pleines de $\hat \Pi$.
Plus généralement, on dira que $B\subset \hat B$ est un {\em
  sous-borélien plein} de $\hat B$ s'il existe un découpage de $\hat
B$ en piles simpliciales $\hat \Pi_j$ tel que $B_j=B\cap \hat \Pi_j$
est un sous-borélien plein de la pile $\hat \Pi_j$ pour tout $j$. On
appellera cela un découpage de $\hat B$ {\em relatif à $B$}, ou bien
un {\em découpage relatif} de la paire $(B,\hat B)$.

\medskip On remarque d'abord le fait élémentaire suivant:

\begin{lem}\label{lem:dec-rel-1}
  Soient $B_1,B_2\subset \hat B$ deux sous-boréliens disjoints de
  $\hat B$. Si $\hat\Pi_i^1$ et $\hat \Pi_j^2$ sont des découpages de
  $\hat B$ relatifs à $B_1$ et $B_2$ respectivement, alors $\hat
  \Pi_{ij}=\hat \Pi_i^1\cap \hat \Pi_j^2$ est un découpage de $\hat B$
  relatif à $B_1\cup B_2$.
\end{lem}

Nous prouvons maintenant un analogue du théorème \ref{thm:dec-pil}
pour des paires de boréliens simpliciaux.

\begin{lem}\label{lem:dec-rel-2}
  Toute paire de boréliens simpliciaux à feuilles compactes $B\subset
  \hat B$ admet un découpage relatif.
\end{lem}
\preuve Soient $T$ et $\hat T$ les espaces de feuilles de $B$ et $\hat
B$ respectivement. On considère l'application $\nu:\hat T\to \N$ qui
assigne à chaque $t\in \hat T$ le nombre de feuilles de $B$ contenues
dans la feuille correspondante $\hat B_t$ de $B$. La mesurabilité de
l'application $\nu$ se montre de façon analogue à celle de
l'application $\sharp$ (voir la proposition \ref{prop:car-mes}). On
peut alors se ramener au cas où la fonction $\nu$ est constante égale
à un nombre entier positif $r(B,\hat B)$ que nous appellerons le {\em
  nombre d'intersection} de la paire.

\medskip\noindent {\bf Première étape.} On suppose $r(B,\hat B)=1$.
Dans ce cas il suffit de montrer que l'application $\hat p \circ i$
dans le diagramme $$
\xymatrix{
  B \ar[r]^i \ar[dr]_{\hat p\circ i} & \hat \Pi \ar[d]^{\hat p}\\
  & \hat \Omega } $$
est semi-simple. Cette BT-application étant
simpliciale, cela découle immédiatement du corollaire
\ref{cor:sim-sem}.

\medskip\noindent {\bf Deuxième étape.} Nous appliquons maintenant une
récurrence sur le nombre d'intersection $r(B,\hat B)$. On suppose que
toute paire de boréliens simpliciaux dont le nombre d'intersection est
plus petit ou égale à $n$ vérifie l'énoncé du lemme, et on considère
une paire $B\subset \hat B$ avec $r(B,\hat B)=n+1$. La projection
naturelle de l'espace de feuilles $T$ de $B$ sur l'espace de feuilles
$\hat T$ de $\hat B$ est borélienne et à fibres finies et admet donc
une section borélienne d'après le théorème \ref{thm:kal}. Celle-ci
découpe $B$ en deux boréliens saturés $B_1$ et $B_2$, le premier
correspondant à l'image de cette section et le deuxième étant son
complémentaire dans $B$. Nous avons alors $r(B_1,\hat B)=1$ et
$r(B_2,\hat B)=n$. Les paires $(B_1,\hat B)$ et $(B_2,\hat B)$
admettent alors des découpages relatifs par hypothèse de récurrence.
Le lemme \ref{lem:dec-rel-1} permet de conclure.  \qed

\subsubsection{Extension d'homéomorphismes locaux.}
Le but de ce paragra\-phe est de montrer le lemme suivant:

\begin{lem}\label{lem:ext-hom}
  Soient $B$ et $\hat B$ deux boréliens simpliciaux par disques, avec
  $B\subset int(\hat B)$. Si $\rho:B\to S$ est un BT-homéomorphisme
  local semi-simple sur une surface fermée $S$, alors il admet une
  BT-extension semi-simple $\hat\rho:\hat B\to S$.
\end{lem}
\preuve Quitte à découper $\hat B$ on peut supposer d'après le lemme
\ref{lem:dec-rel-2} que celui-ci est une pile simpliciale $\hat\Pi$ et
que $B$ se découpe à son tour en exactement $r$ piles
$\Pi_1,\dots,\Pi_r$ qui sont pleines dans $\hat\Pi$. De plus, on peut
supposer que les piles $\Pi_k$ sont contenues dans une des piles du
découpage de $B$ déterminé par l'application semi-simple $\rho$. En
effet, le découpage de $\rho$ sur chaque pile $\Pi_k$ détermine un
découpage de $\hat \Pi$ (voir la propriété (c) en page
\pageref{pil-ple}). L'intersection de ces découpages nous en donne un
dont chaque pile vérifie les conditions requises.

Nous avons donc pour chaque $k=1,\dots,r$ une paire de diagrammes
commutatifs: $$
\xymatrix{
  \Pi_k \ar[r]^{\rho} \ar[d]_{p_k} & S\\
  \Omega_k \ar[ur]_{f_k} } \hspace{1cm} \xymatrix{
  \Pi_k \ar[r]^i \ar[d]_{p_k} & \hat\Pi \ar[d]^{\hat p}\\
  \Omega_k \ar[r]^{g_k} & \hat\Omega } $$
où $p_k$ désigne la
projection de la pile $\Pi_k$ sur sa base $\Omega_k$, $f_k$ est
l'homéomorphisme local associé à l'application simple $\rho|_{\Pi_k}$,
et $g_k$ un plongement simplicial du disque $\Omega_k$ dans
l'intérieur du disque $\bar\Omega$.

Le problème de l'extension semi-simple de $\rho$ est ainsi réduit au
problème de l'extension des homéomorphismes locaux $f_k:\Omega_k\to S$
au disque $\bar\Omega$. Plus précisément, il suffit de trouver un
homéomorphisme local $f:\hat\Omega \to S$ tel que $f_k=f\circ g_k$
pour tout $k$. Pour ce faire on considère le revêtement universel de
$S$: $$
q:\R^2\to S.  $$
On choisit alors des relèvements $\tilde
f_k:\Omega_k\to \R^2$ des $f_k$ relatifs au revêtement $q$, dont les
images sont des disques $D_1,\dots,D_r$ de $\R^2$ deux à deux
disjoints. On considère ensuite un "grand" disque $D$, qui peut être
supposé métrique et centré en $0$, dont l'intérieur contient les
"petits" disques $D_1,\dots,D_r$. On remarque que les surfaces
planaires $$
\hat\Omega-\cup_k\Omega_k\quad,\quad D-\cup_k D_k $$
sont
homéomorphes car elles sont planaires avec le même nombre de
composantes connexes dans le bord. Un homéomorphisme entre ces deux
surfaces permet d'étendre naturellement les plongements $\tilde f_k$
en un plongement $\tilde f:\bar\Omega\to \R^2$. L'homéomorphisme local
$f=q\circ \tilde f$ vérifie alors les conditions requises. Ceci
complète la preuve du lemme.  \qed

\subsubsection{Rétrécissement d'un borélien simplicial.}\label{sec:ret}
On note $sd_n(\K)$ la triangulation de $(X,\mathcal F)$ obtenue au
moyen de la $n$-ième subdivision barycentrique de $\K$. Pour tout
ensemble $A\subset X$, on définit la {\em $n$-étoile} de $A$ comme le
complémentaire $star_n(A)$ des triangles de $sd_n(\K)$ qui ne
rencontrent pas $A$. Il n'est pas difficile de vérifier que si $A$ est
un borélien de $X$, il en est de même pour $star_n(A)$.

Soit $\Omega$ une surface simpliciale pour la triangulation $\K$
contenue dans une feuille de $(X,\mathcal F)$. Pour chaque entier non
négatif $n\in \N$ on définit le {\em $n$-rétrécissement} de $\Omega$
comme étant le fermé $$
ret_n(\Omega)=\Omega-star_n(\partial \Omega).
$$
Il est réunion des triangles de $sd_n(\K)$ contenus dans $\Omega$
et ne rencontrant pas son bord $\partial\Omega$. On remarque que
$ret_n(\Omega)$ est un fermé de $\Omega$ contenu dans l'intérieur de
$ret_{n+1}(\Omega)$ et que $int(\Omega)=\bigcup_n ret_n(\Omega)$. Pour
$n\geq 1$, le $n$-rétrécissement d'un disque simplicial pour $\K$ est
un disque simplicial pour $sd_n(\K)$.

Soit $B$ un borélien simplicial de $(X,\mathcal F)$. On définit le
{\em $n$-rétrécisse\-ment} de $B$ comme étant l'ensemble $ret_n(B)$
obtenu par $n$-rétrécisse\-ment de chaque feuille de $B$. Il est très
facile de voir (en tenant compte du théorème \ref{thm:dec-pil}) que si
$B$ est à feuilles compactes alors $ret_n(B)$ est aussi un borélien
(simplicial pour la triangulation $sd_n(\K)$ et à feuilles compactes).
En effet, le $n$-rétrécissement d'une pile simpliciale est clairement
une pile simpliciale.

\subsubsection{Preuve de la proposition \ref{prop:lam-pla-rev}.}
La preuve consiste en une application répétée du lemme
\ref{lem:ext-hom}.  Reprenons la filtration par disques $B_n$ de la
lamination $(X,\mathcal F)$.  Pour chaque $n$ le borélien $B_n$ est
contenu dans $B_{n+1}$, mais il n'est pas contenu dans l'intérieur de
$B_{n+1}$ La paire $(B_n,B_{n+1})$ ne vérifie donc pas les hypothèses
du lemme \ref{lem:ext-hom}. Pour pallier à cette difficulté il suffit
de rétrécir convenablement les $B_n$. On définit $$
\hat
B_n=ret_n(B_n)\quad,\quad n\in \N $$
Les boréliens $\hat B_n$
vérifient les conditions d'une filtration compacte forte à l'exception
du fait qu'il ne sont pas tous simpliciaux relativement à la même
triangulation. Pour tout $n$ la paire $(\hat B_n,\hat B_{n+1})$
vérifie tout de même les conditions du lemme \ref{lem:ext-hom}
relativement à la triangulation $sd_{n+1}(\K)$.

On se fixe un découpage en piles de $B_1$ et on se donne un premier
BT-homéomorphisme local $\rho_1:B_1\to S$, défini sur chaque pile
$\Pi$ de $B_1$ par le diagramme suivant: $$
\xymatrix{
  \Pi \ar[r]^{\rho_1} \ar[d]_p & S\\
  \Omega \ar[ur]_{f} } $$
où $p$ est la projection de $\Pi$ sur sa
base $\Omega$, et $f:\Omega\to S$ un homéomorphisme local quelconque.
Le BT-homéomorphisme local $\rho_1$ est alors semi-simple par
construction. On applique le lemme \ref{lem:ext-hom} à chaque étape
pour obtenir une suite $\rho_n:\Pi_n\to S$ de BT-homéomor\-phismes
locaux semi-simples, chacun étendant le précédent. La limite
$\rho=\lim \rho_n$ est alors un BT-homéomorphisme local globalement
défini.

Rien n'a été pourtant dit jusqu'à présent qui permette d'assurer que
$\rho$ est un BT-revêtement. Pour ce faire il suffit de bien choisir,
dans la preuve de \ref{lem:ext-hom}, les disques métriques $D$, de
sorte que ceux correspondants à l'étape $n$ aient un rayon supérieur
ou égal à $n$. De cette façon on assure la complétude du
BT-homéomorphisme local $\rho$ qui devient donc un BT-revêtement. La
proposition est démontrée.

\begin{rem}
  Si $(X,\mathcal F)$ est une lamination par plans hyperfinie de
  classe $C^p$, avec $0\leq p \leq \infty$, on peut construire le
  BT-revêtement $\rho$ de classe $C^p$. Les détails sont laissés au
  lecteur.
\end{rem}

\subsection{Laminations par cylindres}
Une lamination borélienne $(X,\mathcal F)$ est dite {\em par
  cylindres} si toutes ses feuilles sont homéomorphes à $\mathbb
S^1\times \R$.

\begin{thm}\label{cor:rev}
  Soit $(X,\mathcal F)$ une lamination borélienne par cylindres. Alors
  $(X,\mathcal F)$ est un BT-revêtement du tore $\mathbb T^2$.
\end{thm}

On remarque d'abord que {\em toute lamination borélienne par cylindres
  $(X,\mathcal F)$ est hyperfinie}, car l'extension compacte d'un flot
\cite{Gh1}.  Plus précisément, il existe une lamination borélienne
$(X_1,\mathcal F_1)$ de dimension $1$ et une BT-application $X\to X_1$
à fibres compactes. La lamination $(X_1,\mathcal F_1)$ étant
hyperfinie, il en a de même pour $(X,\mathcal F)$. Ce résultat a été
prouvé par Ghys dans \cite{Gh1} pour le cas d'une laminations
topologique munie d'une mesure harmonique, et il a été adapté au cas
d'une mesure quasi-invariante quelconque par Blanc dans \cite{Bla}.
Enfin il peut être facilement adapté au cadre purement borélien en
utilisant des triangulations.

La preuve du théorème \ref{cor:rev} suit la même démarche que celle de
du théorème \ref{thm:rev}, mais en remplaçant les filtrations
compactes par disques par des filtrations compactes dont les feuilles
sont des cylindres compacts $\mathbb S^1\times [0,1]$. Les détails
seront laissés au lecteur.

\begin{rem} {\em Le type d'une lamination par plans}.  Une autre
  question se pose naturellement à la vue du théorème \ref{thm:rev}:
  est-ce que toute lamination par plans (non hyperfinie) est un
  BT-revêtement d'une surface compacte? Le théorème \ref{thm:rev} ne
  répond pas à cette question car il est possible de construire des
  laminations par plans non hyperfinies (voir l'exemple du
  \S\ref{ex:non-hyp}). En général, la suspension d'une action
  borélienne libre et préservant une mesure de Borel finie du groupe
  fondamental de toute surface de genre $\geq 2$ est une lamination
  par plans non hyperfinie. Or tous ces exemples sont des
  BT-revêtements de la surface de genre $2$, mais le genre de celle ci
  ne peut pas être réduit d'après le théorème \ref{thm:rev}. Nous
  appellerons {\em type} d'une lamination orientable $(X,\mathcal F)$
  le plus petit nombre réel $n$ tel que $(X,\mathcal F)$ est un
  BT-revêtement de la surface compacte orientable de genre $n$. Si se
  nombre n'existe pas, on dira $(X,\mathcal F)$ de {\em type infini}.
  On remarquera que, toute surface orientable étant le revêtement
  d'une surface de genre $0$, $1$ ou $2$, le type d'une lamination par
  plans ne peut être que $0$, $1$, $2$ ou $\infty$. D'après le
  théorème \ref{thm:rev}, {\em une lamination par plans est hyperfinie
    si et seulement elle est de type $1$}. Nous ne connaissons pas
  d'exemples de laminations par plans de type $\infty$.
\end{rem}

\section{Laminations paraboliques}\label{sec:par}
\subsection{Le théorème d'uniformisation}
Nous rappelons le théorème classique suivant:

\begin{thm}[Koebe-Poincaré]\label{koebe}
  Toute surface de Riemann simplement connexe est conformément
  équivalente à l'une de ces trois surfaces: la sphère de Riemann, le
  plan euclidien ou le disque de Poincaré.
\end{thm}

En particulier pour toute surface $S$ munie d'une métrique de Riemann
$g$, il existe une fonction $u:S\to \R$ telle que la métrique
$\exp(u)g$ est à courbure constante égale à $+1$, $0$ ou $-1$. Selon
la valeur de cette constante on dira que $S$ est {\em elliptique},
{\em parabolique} ou {\em hyperbolique} respectivement.

\begin{ex}
  Une lamination par surfaces $(X,\mathcal F)$ est dite {\em
    holomorphe} si elle est définie par un atlas dont les changements
  de cartes sont holomorphes. En particulier les feuilles de
  $(X,\mathcal F)$ sont munies d'une structure conforme. Toute
  lamination holomorphe $(X,\mathcal F)$ peut être munie d'une
  métrique de Riemann conforme. Réciproquement, toute métrique de
  Riemann détermine un atlas holomorphe donné par des systèmes locaux
  de coordonnées isothermes. On dira qu'une lamination holomorphe
  $(X,\mathcal F)$ est {\em parabolique} si toutes ses feuilles sont
  paraboliques. Rappelons que les feuilles d'une telle lamination sont
  des plans, des cylindres ou des tores.
\end{ex}

Le lemme suivant est une généralisation du théorème de Koebe-Poincaré
aux familles de structures conformes paraboliques. Un résul\-tat
analogue se vérifie pour les familles de métriques hyperboliques et
elliptiques, mais nous n'en ferons pas l'usage. Il peut être aussi vu
comme une version mesurable du théorème d'Ahlfors-Bers \cite{AB}.

\begin{lem}\label{lem:par-tri}
  Soit $(X,\mathcal F)$ une lamination parabolique telle qu'il existe
  un BT-isomorphisme $X\simeq \R^2\times T$ entre $X$ et un prisme de
  base $\R^2$.  Alors il existe une BT-application: $$
  f:X\to \C $$
  qui est une équivalence conforme en restriction à chaque feuille de
  $(X,\mathcal F)$.
\end{lem}

\preuve La donnée de $f$ équivaut à celle d'une famille boré\-lienne
d'équivalences conformes $f_t:L_t\to \C$ où $L_t$ est la feuille de
$(X,\mathcal F)$ correspondant au point $t\in T$.

On utilise les cartes locales d'un atlas feuilleté holomorphe de
$(X,\mathcal F)$ pour fixer de façon borélienne un paramétre complexe
$z_t$ au voisinage de $0$ dans chaque feuille $L_t$. D'après le
théorème de Koebe-Poincaré, il existe une bijection holomorphe non
constante $f_t:L_t\to \C$ pour tout $t\in T$. Cette fonction est
donnée par $f_t=M_t\circ g_t$, où $g_t:L_t\to \C\cup\{\infty\}$ est
une fonction méromorphe non constante et $M_t$ une transformation de
M\"obius (cf. \cite{FK}, \S IV.4.8). La fonction $g_t$ est définie
sous la forme: $$
g_t=u_t+\sqrt{-1}v_t $$
où $u_t$ est une fonction
harmonique dans $L_t-\{0\}$ limite d'une suite $u_t^n$, chacune
solution du problème de Dirichlet dans $M_t-\{|z_t|<1/n\}$ donné par
la condition de contour $u_t^n(z)=\mathrm{Re}(1/z_t)$. La fonction
$v_t$ est limite de solutions du problème de Dirichlet sur
$L_t-\{|z_t|<1/n\}$ avec $\mathrm{Re}(\sqrt{-1}/z_t)$ comme condition
de contour (voir \cite{FK}, \S IV.4 pour les détails).

Nous remarquons que la méthode de Perron (voir par exemple \cite{AS},
\S III.3) entraîne la dépendance me\-surable des solutions du problème
de Dirichlet relativement aux conditions de contour. En effet, si $U$
est une région régulière de $M$ et $c$ est une fonction complexe
continue définie sur $\partial U$, alors il existe une seule extension
de $c$ à $\overline U$ qui soit harmonique sur $U$. Cette fonction est
donnée par la formule universelle: $$
u(z)=\sup_{v\in \mathcal V(c)}
v(z) $$
où $\mathcal V(c)$ est la famille des fonctions
sous-harmoniques sur $U$ telles que: $$
\varlimsup_{z\in \zeta}
v(z)\leq c(\zeta) $$
pour tout $\zeta\in \partial U$. Ceci montre la
mesurabilité des familles $u_t$ et $v_t$, et donc celle de $f_t$.\qed

\bigskip
\subsection{Groupoïde fondamental}
Soit $(X,\mathcal F)$ une lamination borélienne de dimension $n$. Un
chemin de $(X,\mathcal F)$ est une application continue $c:[0,1]\to L$
où $L$ est une feuille de $(X,\mathcal F)$. On note $\Pi(X,\mathcal
F)$ l'espace formé par les classes d'homotopie de chemins de
$(X,\mathcal F)$. Rappelons que les homotopies considérées sont
supposées préserver l'origine et l'extrémité des che\-mins. Les
applications qui associent à tout chemin tangent son origine et son
extrémité déterminent donc des applications: $$
\xymatrix{
  \Pi(X,\mathcal F)\ar@/^/[rr]^\alpha \ar@/_/[rr]_\beta & & X } $$
bien définies au niveau des classes d'homotopie. L'espace
$\Pi(X,\mathcal F)$ est naturellement muni d'une structure de
groupoïde dont $\alpha$ et $\beta$ sont les applications source et but
respectivement, et dont la loi de composition est donné par passage
aux classes d'homotopie de la somme usuelle de chemins. On appelle
$\Pi(X,\mathcal F)$ le {\em groupoïde fondamental} de la lamination
$(X,\mathcal F)$.

Pour tout $x\in X$, on notera $\tilde L_x$ et $\tilde L^x$ la
$\alpha$-fibre et la $\beta$-fibre de $x$. Tout élément $\gamma\in
\Pi(X,\mathcal F)$ avec $x=\alpha(\gamma)$ et $y=\beta(\gamma)$ agit à
gauche sur les points de $\tilde L_y$ les envoyant sur $\tilde L_x$.
De façon analogue, il agit à droite sur les points de $\tilde L^x$ les
envoyant sur $\tilde L^y$. On notera ces applications: $$
\tau_\gamma:\tilde L_y\to \tilde L_x \quad,\quad \rho_\gamma:\tilde
L^x\to \tilde L^y $$

On peut définir de façon analogue le groupoïde fondamental $\Pi(X)$
d'un BT-espace standard quelconque $X$ et le munir d'une
BT-struc\-ture dont les feuilles sont les groupoïdes fondamentaux des
feuilles de $X$. Dans le cas d'une lamination borélienne, nous
définirons explicitement la BT-structure du groupoïde fondamental par
le biais des cartes d'un atlas feuilleté borélien.

\medskip
\subsubsection{La structure BT de $\Pi(X,\mathcal F)$.}
Par analogie avec la structure topologique ou différentiable des
groupoïdes des feuilletages (cf. \cite{AH}) nous allons définir une
lamination borélienne sur l'espace $\Pi(X,\mathcal F)$ dont les
feuilles sont des variétés de dimension $2n$. Cette lamination est
l'analogue borélien des laminations induites $\alpha^*(\mathcal
F)=\beta^*(\mathcal F)$ dans le cas différentiable. Les applications
source et but seront des BT-submersions dont les fibres définiront des
laminations de dimension $n$ sur $\Pi(X,\mathcal F)$ symétriques
relativement à la structure de groupoïde.

Soit $\mathbb B^n$ la boule ouverte unité de $\R^n$ et $T$ un espace
de Borel standard. Soit $\mathbb B^n\times T\times \mathbb B^n$ un
prisme de verticale $T$.  Une BT-application: $$
\mathfrak
c:[0,1]\times \mathbb B^n\times T\times \mathbb B^n \to X $$
est dite
un {\em tube de chemins} de $(X,\mathcal F)$ si les applications: $$
\mathfrak c_0(\cdot,\cdot,y):\mathbb B^n\times T\to X \quad,\quad
\mathfrak c_1(x,\cdot,\cdot):T\times \mathbb B^n\to X $$
sont des
BT-plongements quels que soient $x,y\in \mathbb B^n$. Les deux
applications ci-dessus représentent des paramétrages locaux de la
lamination $(X,\mathcal F)$. Les applications $\mathfrak
c(\cdot,x,t,y):[0,1]\to X$ sont des chemins de $(X,\mathcal F)$ qui
relient les points $\mathfrak c_0(x,t,y)$ et $\mathfrak c_1(x,t,y)$.
Par passage aux classes d'homotopie de chemins nous obtenons une
application injective: $$
\hat {\mathfrak c}:\mathbb B^n\times T\times
\mathbb B^n\to \Pi(X,\mathcal F) $$
On considère une famille
dénombrable de tubes de chemins $\mathfrak c_i$ et on note $U_i$
l'image dans $\Pi(X,\mathcal F)$ de chaque $\hat {\mathfrak c}_i$.  Si
la famille $U_i$ recouvre $\Pi(X,\mathcal F)$ alors les applications
$$
\tilde\varphi_i=\hat {\mathfrak c}_i^{-1}:U_i\to \mathbb B^n\times
T\times \mathbb B^n $$
sont les cartes d'un atlas feuilleté borélien
sur $\Pi(X,\mathcal F)$. On notera $\tilde {\mathcal F}$ la lamination
borélienne définie par cet atlas.  Ceci munit l'espace $\Pi(X,\mathcal
F)$ d'une BT-structure, de sorte que l'inclusion canonique: $$
X\to
\Pi(X,\mathcal F) $$
qui consiste à envoyer le point $x\in X$ sur la
classe d'homotopie du lacet constant en $x$, est en fait un plongement
de BT-espaces.

On vérifie aisément que la lamination $(\Pi(X,\mathcal F),\tilde
{\mathcal F})$ peut être donné par un atlas feuilleté de la même
classe que $(X,\mathcal F)$.  Les applications $$
\xymatrix{
  \Pi(X,\mathcal F)\ar@/^/[rr]^\alpha \ar@/_/[rr]_\beta & & X } $$
sont des BT-applications dont la restriction à chaque feuille de
$\tilde {\mathcal F}$ est une submersion sur une feuille de
$(X,\mathcal F)$. Remarquons d'ailleurs que les feuilles de cette
lamination $(\Pi(X,\mathcal F),\tilde {\mathcal F})$ ne sont que les
groupoïdes fondamentaux des feuilles de $(X,\mathcal F)$. Les fibres
de $\alpha$ (resp.  $\beta$) définissent sur $\Pi(X,\mathcal F)$ une
autre lamination plus fine de même classe et dimension que $\mathcal
F$ que nous noterons $\mathcal F_\alpha$ (resp. $\mathcal F_\beta$).
Remarquons que les translations à gauche et à droite: $$
\tau_\gamma:\tilde L_y\to \tilde L_x \quad,\quad \rho_\gamma:\tilde
L^x\to \tilde L^y $$
sont des difféomorphismes de classe maximale
entre les feuilles de $\mathcal F_\alpha$ et $\mathcal F_\beta$
respectivement.

\subsection{Moyennabilité}
Une {\em moyenne} sur un espace de Borel standard mesuré $(T,\lambda)$
est par définition un état dans $L^\infty(T,\lambda)$. Autrement dit,
une moyenne est une fonctionnelle continue $m:L^\infty(T,\lambda)\to
\R$ positive sur les fonctions positives, et telle que $m({\bf 1})=1$.

\begin{defn}
  Une lamination borélienne différentiable $(X,\mathcal F)$ est dite
  {\em moyennable} si elle admet une {\em moyenne} relative à sa
  classe de Lebesgue, i.e. une application linéaire positive: $$
  m:L^\infty(X,\mathcal F)\to L^\infty(X,\mathcal F) $$
  qui à tout
  élément $f\in L^\infty(X,\mathcal F)$ associe une fonction $m(f)\in
  L^\infty(X,\mathcal F)$ basique, i.e. constante le long des
  feuilles, et telle que $m(\bf 1)=\bf 1$. Ici $L^\infty(X,\mathcal
  F)$ désigne l'espace vectoriel des fonctions boréliennes $f:X\to
  \mathcal \R$ qui sont $\lambda$-essentiellement bornées sur chaque
  feuille, où $\lambda$ est la classe de Lebesgue de $(X,\mathcal F)$.
\end{defn}

Pour des raisons pratiques, nous allons donner une définition
équi\-valente à celle de moyenne: celle de système borélien de
moyennes. Un {\em système borélien de moyennes} sur $(X,\mathcal F)$
est une fonction $m$ qui assigne à chaque $x\in X$ une moyenne: $$
m_x\in (L^\infty(\tilde L_x,\lambda_x))^* $$
sur la feuille
correspondante vérifiant les conditions suivantes:
\begin{enumerate}
\item $m_x=m_y$ pour $x$ et $y$ dans la même feuille de $X$;

\item Pour toute fonction borélienne bornée $f:X\to \R$ la fonction
  $x\mapsto m_x(f)$ est borélienne.
\end{enumerate}

Un système borélien de moyennes $m_x$ sur $X$ définit sur celui-ci une
moyenne $m$ en posant $m(f)(x)=m_x(f)$. De façon analogue, une moyenne
définit un système borélien de moyennes par restriction aux feuilles.
Les deux notions sont donc équivalentes.

\begin{defn}\label{def:gro-moy}
  Soit $(X,\mathcal F)$ une lamination borélienne différentiable et
  soit $\Pi(X,\mathcal F)$ son groupoïde fondamental. On dira que
  $\Pi(X,\mathcal F)$ est {\em moyennable à gauche} s'il existe un
  système borélien de moyennes $m$ sur la lamination $\tilde{\mathcal
    F_\alpha}$ tel que
  \begin{itemize}
  \item Pour tout $\gamma\in \Pi(X,\mathcal F)$ avec
    $x=\alpha(\gamma)$ et $y=\beta(\gamma)$, et toute fonction
    borélienne $f\in L^\infty(\tilde L_x,\mu)$, on a: $$
    \tilde
    m_x(f)=\tilde m_y(\tau^*_\gamma(f)) $$
  \end{itemize}
  On définit de façon analogue la moyennabilité à droite.
\end{defn}

\begin{rem}
  La définition de moyennabilité apparaît dans le cas des groupes
  localement compacts (cf. \cite{Gre}). Si $G$ est un groupe
  localement compact et $\nu$ sa mesure de Haar, alors on définit une
  moyenne invariante à gauche (resp. à droite) sur $G$ comme étant un
  état dans $L^\infty(G,\nu)$ invariant par l'action de $G$ sur
  l'espace $(L^\infty(G,\nu))^*$ induite par l'action à gauche (resp.
  à droite) de $G$ sur lui même. Cette définition a été généralisée
  par Anantharaman-Delaroche et Renault \cite{AR} au cas des
  groupoïdes boréliens, eux aussi munis d'un "système de mesures de
  Haar" et d'une mesure transverse quasi-invariante. Notre définition
  est plus forte la définition 3.2.1 de \cite{AR}. Notre définition
  \ref{def:gro-moy} est en effet celle de \cite{AR}, mais à un détail
  près: la mesurabilité du système de moyennes est considérée ici de
  manière globale et non plus relativement à une mesure
  quasi-invariante. La moyennabilité d'un groupoïde au sens de notre
  définition entraîne la moyennabilité du groupoïde au sens de
  \cite{AR} quelle que soit la mesure quasi-invariante choisie, ce que
  Anantharaman-Delaroche et Renault appellent la {\em measurewise
    amenability} du groupoïde.
\end{rem}

\begin{ex}\label{ex:BT-act}
  Si $(X,\mathcal F)$ est définie par une action $\phi$ d'un groupe de
  Lie moyennable $G$, alors $(X,\mathcal F)$ est moyennable. En effet,
  si $m\in (L^\infty(G,\nu))^*$ est une moyenne invariante à gauche
  sur $G$, alors nous avons sur $X$ un système borélien de moyennes
  $\phi_*(m)$ défini par $$
  \phi_*(m)_x(f)=m(\phi(\cdot,x)^*f)\quad,\quad f\in
  L^\infty(L_x,\lambda_x).  $$
\end{ex}

Cet exemple permet d'illustrer le fait suivant:

\begin{prop}
  Soit $(X,\mathcal F)$ une lamination différentiable. Si son
  groupoïde fondamental $\Pi(X,\mathcal F)$ est moyennable à gauche ou
  à droite, alors $(X,\mathcal F)$ est moyennable.
\end{prop}
\preuve En effet, si $\tilde m$ est un système de moyennes invariant à
gauche sur $\Pi(X,\mathcal F)$, on définit $m_x(f)=\tilde
m_x(\beta^*_x)$, où $\beta^*_x:\tilde L_x\to L_x$ est le revêtement
conforme de la feuille $L_x$ défini par l'application but $\beta$. Le
même raisonnement s'applique au cas de la moyennabilité à droite.
\qed

\medskip Remarquons toutefois que la réciproque n'est pas vraie en
général, c'est-à-dire qu'il existe des laminations moyennables dont le
groupoïde fondamental est non moyennable.

\subsection{Moyennabilité des laminations paraboliques}
On retourne ici au cas des laminations par surfaces, et plus
particulièrement au cas des laminations holomorphes. Ce sont par
définition des laminations différentiables, de sorte que sa classe de
Lebesgue est bien définie.

\medskip
\begin{thm}\label{thm:par-moy}
  Si $(X,\mathcal F)$ est une lamination holomorphe parabolique, alors
  son groupoïde fondamental $\Pi(X,\mathcal F)$ est moyennable à
  gauche et à droite.
\end{thm}

\preuve Les feuilles de la lamination $\mathcal F_\alpha$, définie sur
$\Pi(X,\mathcal F)$, sont toutes des plans car elles peuvent être
identifiés au revêtement universel des feuilles de $\mathcal F$. Elle
est d'ailleurs clairement hyperfinie car elle admet un domaine
fondamental, i.e. un borélien de $\Pi(X,\mathcal F)$ qui rencontre
chaque feuille en un seul point. Elle est donc, d'après le théorème
\ref{thm:rev}, un BT-revêtement du tore $\mathbb T^2$.

L'action naturelle de $\R^2$ sur celui-ci se relève en une BT-action
localement libre de $\R^2$ sur $\Pi(X,\mathcal F)$ $$
\phi:\R^2\times
\Pi(X,\mathcal F) \to \Pi(X,\mathcal F) $$
définissant la lamination
$\mathcal F_\alpha$. Sa restriction à $X$ donne un BT-isomor\-phisme:
$$
\phi|_{\R^2\times X}:\R^2\times X \to \Pi(X,\mathcal F).  $$
Par
conséquent le BT-espace $\Pi(X,\mathcal F)$ est un prisme et
$(\Pi(X,\mathcal F),\mathcal F_\alpha)$ est vérifie les conditions du
lemme \ref{lem:par-tri}. Il existe une BT-application $$
f:\Pi(X,\mathcal F)\to \C $$
qui est une équivalence conforme en
restriction à chaque feuille de $\mathcal F_\alpha$.

Le groupe $G$ des similitudes du plan euclidien $\C$ est un produit
semi-direct de groupes abéliens, donc moyennable. Son action sur $\C$
étant transitive et préservant la classe de la mesure de Lebesgue
$\lambda$, il existe une moyenne $\tilde m\in
(L^\infty(\C,\lambda))^*_1$ invariante par $G$. Le relèvement de la
moyenne $\tilde m$ par le biais de l'application $f$ nous donne un
système borélien de moyennes $\tilde m_x\in(L^\infty(\tilde
L_x,\lambda))^* $, $x\in X$, sur la lamination $\mathcal F_\alpha$.
Mais rappelons que, pour tout $\gamma\in\Pi(X,\mathcal F)$, la
translation à gauche $\tau_\gamma:\tilde L_y \to \tilde L_x$ est une
équivalence conforme. L'invariance de $\tilde m$ par l'action du
groupe conforme $G$ entraîne alors celle de $\tilde m_x$ par les
translations $\tau_\gamma$.  En d'autres mots, nous avons une moyenne
invariante à gauche sur le groupoïde $\Pi(X,\mathcal F)$. On raisonne
de façon analogue pour montrer la moyennabilité à droite. Le théorème
est démontré.  \qed

\subsection{Revêtements du tore, métriques plates et actions localement libres
de $\R^2$}\label{sec:rev-tor}
Le but de ce paragraphe est de montrer le théorème \ref{thm:eq-mes}
énoncé dans l'introduction. On supposera que $(X,\mathcal F)$ est une
lamination borélienne orientable de classe $C^1$ munie d'une mesure
transverse quasi-invariante $\mu$.

Supposons que $(X,\mathcal F)$ est hyperfinie par plans, tores et
cylindres mod$(\mu)$. Ceci signifie qu'il existe un borélien saturé de
mesure totale $X'\subset X$ tel que la lamination induite
$(X',\mathcal F')$ est hyperfinie par plans, tores et cylindres. On
peut découper maintenant $X'$ en trois boréliens saturés $X'_1$,
$X'_2$ et $X'_3$ réunion respectivement des tores, cylindres et plans
de la lamination $\mathcal F'$ (voir par exemple \cite{Gh1}).  Les
laminations induites par $\mathcal F'$ sur $X'_2$ et $X'_3$ sont des
BT-revêtements de $\mathbb T^2$ d'après les théorèmes \ref{thm:rev} et
\ref{cor:rev}. Enfin la lamination par tores induite par $\mathcal F'$
sur $X'_1$ est un BT-revêtement car $X'_1$ admet un découpage en piles
(cf. théorème \ref{thm:dec-pil}). L'implication (1)$\Rightarrow$(2)
est démontrée.

Un BT-revêtement de $\mathbb T^2$ est défini par une BT-action
localement libre de $\R^2$ relèvement de l'action naturelle de
celui-ci sur $\mathbb T^2$. Une telle action détermine de manière
évidente une métrique plate complète sur les feuilles, qui est
parabolique d'après le théorème d'uniformisation de Koebe-Poincaré.
Ceci montre les implications
(2)$\Rightarrow$(3)$\Rightarrow$(4)$\Rightarrow$(5).

Enfin si une lamination $(X,\mathcal F)$ admet une métrique de Riemann
parabolique mod$(\mu)$ alors elle est moyennable mod$(\mu)$ d'après le
théo\-rème \ref{thm:par-moy}. En particulier elle est $\mu$-moyennable
et donc hyperfinie mod$(\mu)$ d'après \cite{CFW}. Ses feuilles sont
des plans, tores ou cylindres car ce sont les seules trois surfaces
admettant une métrique parabolique.

\bigskip
\section{Existence de triangulations}\label{sec:exi-tri}
Nous esquissons ici la preuve du fait que toute lamination borélienne
$(X,\mathcal F)$ de dimension $2$ admet une triangulation. Rappelons
qu'une triangulation est donnée par une famille dénom\-brable de piles
de triangles $$
\pi_i:\Delta\times T_i \simeq \Sigma_i $$
de
$(X,\mathcal F)$. Chaque triangle $\sigma$ est muni d'un paramétrage
canonique $\pi_\sigma:\Delta\to \sigma$. On demande les deux
conditions suivantes:
\begin{itemize}
\item Si $\sigma$ et $\tau$ sont deux triangles différents
  d'intersection non vide, alors ils ont soit une arête, soit un
  sommet en commun.

\item Si $\sigma\cap \tau$ se rencontrent le long d'une arête, alors
  le changement de paramétrage $\pi_\sigma^{-1}\circ \pi_\tau$ est une
  bijection affine entre deux arêtes de $\Delta$.
\end{itemize}

On supposera que la lamination $(X,\mathcal F)$ est de classe $C^1$,
i.e. elle est donnée par un atlas feuilleté borélien: $$
\varphi_i:U_i\to \mathbb B\times T_i $$
dont les changement de plaques
$f_{ij}(\cdot,t)$ sont des difféomorphismes locaux de classe $C^1$.
Ici $\mathbb B$ désigne la boule ouverte unité de $\R^2$.

\subsection{Existence de métriques complètes}
Quitte à prendre un raffinement convenable des $U_i$ on peut supposer
que les plaques de ceux-ci sont relativement compactes dans les
feuilles de $(X,\mathcal F)$. On définit alors une métrique
$g_i=\varphi_i^*|dz|$ sur les plaques de chaque $U_i$. Une partition
de l'unité permet de recoller ces métriques pour obtenir une métrique
globalement définie $g$ sur les feuilles de $\mathcal F$. Puisque les
plaques des $U_i$ sont relativement compactes, la métrique $g$ est
complète le long de chaque feuille. Nous avons donc le résultat
suivant.

\begin{lem}
  Toute lamination borélienne de classe $C^1$ admet une métrique
  complète.
\end{lem}

\subsection{Existence d'un atlas adapté}
Soit $g$ une métrique complète sur une lamination borélienne
$(X,\mathcal F)$.  Un atlas feuilleté borélien de $(X,\mathcal F)$ est
dit {\em convexe} si toutes ses plaques sont convexes.

\begin{lem}\label{lem:atl-con}
  Toute lamination borélienne admet un atlas convexe.
\end{lem}
\preuve Pour tout $x\in X$ et tout $\epsilon\in \R^+$, on note
$B_\epsilon(x)$ la boule de rayon $\epsilon$ centrée en $x$. Il est
bien connu que pour tout $x$ il existe un rayon $\epsilon(x)>0$ tel
que $B_\delta(x)$ est un voisinage convexe de $x$ pour tout $\delta
\leq \epsilon(x)$. Nous laissons au lecteur la tâche très simple de
vérifier que l'on peut supposer que la fonction $\epsilon(x)$ est
borélienne.

Soit $T$ transversale borélienne de $(X,\mathcal F)$. Il est possible
de découper $T$ en une famille $T_i$ de sorte que, pour tout $i$ et
toute paire $x,y\in T_i$, les boules $B_{\epsilon(x)}(x)$ et
$B_{\epsilon(y)}(y)$ sont disjointes. On posera: $$
U_i=\bigcup_{x\in
  T_i} B_{\epsilon(x)}(x) $$
le borélien réunion des $\epsilon$-boules
centrés en les points de $T_i$. On définit alors un paramétrage $$
\pi_i:\mathbb B\times T_i \to U_i $$
en normalisant l'application
exponentielle relativement à la fonction $\epsilon$. Si l'on choisit
la transversale $T$ suffisamment grande pour que la famille de
boréliens $U_i$ soit un recouvrement de $X$, alors les paires
$(U_i,\pi_i^{-1})$ constituent par construction un atlas feuilleté
borélien convexe de $(X,\mathcal F)$.  \qed

Soit $\mathcal A$ un atlas feuilleté de $(X,\mathcal F)$. On appelle
l'{\em étoile} d'un point $x\in X$ relativement à $\mathcal A$
l'ouvert $st_{\mathcal A}(x)$ réunion des plaques de $\mathcal A$ qui
contiennent $x$. Un atlas sera dit {\em adapté} à une métrique $g$ si
$st_{\mathcal A}(x)$ est un convexe de $g$ pour tout $x\in X$. Nous
avons le résultat suivant, dont la preuve est analogue à celle de
\ref{lem:atl-con}.

\begin{lem}\label{lem:atl-ada}
  Toute lamination borélienne admet un atlas adapté.
\end{lem}

\subsection{Construction d'une triangulation}
\begin{figure}[H]
  \begin{center}
    \input{triangulations.pstex_t}
  \end{center}
  \caption{\label{fig:tri-lam} Triangulation des laminations}
\end{figure}
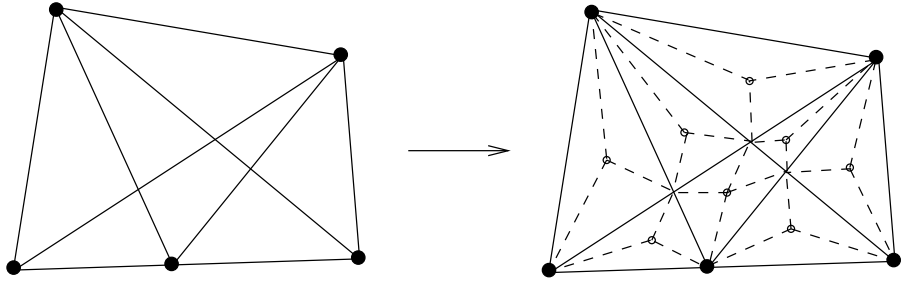

Soit $(X,\mathcal F)$ une lamination borélienne et $g$ une métrique
complète.  On se fixe un atlas $(U_i,\varphi_i)$ de $(X,\mathcal F)$
adapté à la métrique $g$. On appelle {\em centre} d'une plaque $U_t$
de $U_i$ le point $b_t=\varphi_i^{-1}(0,t)$, où $0$ est le centre de
la boule $\mathbb B$. Pour chaque paire de plaques $U_t$ et $U_s$ qui
se rencontrent, on trace la seule géodésique qui rejoint les centres
$b_t$ et $b_s$. Nous obtenons localement un dessin illustré par la
figuré \ref{fig:tri-lam}.

\bigskip Nous avons ainsi une décomposition des feuilles de
$(X,\mathcal F)$ en disques fermés, délimités par des morceaux de
géodésique. Nous fixons dans chaque disque son centre de gravité
relatif au volume de la métrique $g$, et nous traçons les géodésiques
qui rejoignent le barycentre de chaque disque avec les centres des
plaques contenus dans ce disque. Le résultat est une famille de
triangulations $\{\K_L~|~L\in \mathcal F\}$ des feuilles de
$(X,\mathcal F)$.  Chaque triangle peut être identifié à son centre de
gravité. L'espace $\K$ de tous les triangles devient ainsi une
transversale borélienne de $(X,\mathcal F)$.

On décompose $\K$ en une famille dénombrable de boréliens $T_i$ de
sorte que deux triangles différents d'un même $T_i$ ne se coupent pas.
On note $$
\Sigma_i=\bigcup_{\sigma\in T_i}\sigma $$
le borélien
réunion de tous les triangles de $T_i$ et on construit des
paramétrages $$
\pi_i:\Delta\times T_i\to \Sigma_i $$
qui envoient
homothétiquement les arêtes de $\Delta$ sur celles de $\K$. Le
coefficient de cette homothétie est unique et égal à la longueur de
l'arête correspondante. Les changements de paramétrage sont affines
par construction. La famille de piles $\Sigma_i$ est la famille
génératrice d'une triangulation de $(X,\mathcal F)$.

\bigskip
\begin{center}
  \begin{small}{\sc Miguel Bermúdez, Gilbert Hector}\\
    Institut Girard Desargues\\
    Université Claude Bernard Lyon 1\\
    Bâtiment Braconnier\\
    21 avenue Claude Bernard\\
    69622 Villeurbanne cedex\\
    France\\
    E-mail: bermudez@igd.univ-lyon1.fr,
    hector@igd.univ-lyon1.fr\end{small}
\end{center}

\end{document}

%% file: triangulations.pstex_t
\begin{picture}(0,0)%
\includegraphics{triangulations.pstex}%
\end{picture}%
\setlength{\unitlength}{4144sp}%
\begingroup\makeatletter\ifx\SetFigFont\undefined%
\gdef\SetFigFont#1#2#3#4#5{%
  \reset@font\fontsize{#1}{#2pt}%
  \fontfamily{#3}\fontseries{#4}\fontshape{#5}%
  \selectfont}%
\fi\endgroup%
\begin{picture}(5361,1654)(472,-2045)
\end{picture}%